\documentclass[a4paper,11pt]{article}%
\usepackage{amsmath}
\usepackage{graphicx}
\usepackage{amsfonts}
\usepackage{amssymb}%
\begin{document}

\title{Schubert calculus and cohomology of Lie groups Part I. 1--connected Lie groups}
\author{Haibao Duan\thanks{The first author is supported by 973 Program 2011CB302400
and NSFC 11131008.}\ \\Institute of Mathematics, Chinese Academy of Sciences\\dhb@math.ac.cn
\and Xuezhi Zhao\thanks{The second author is supported by NSFC 10931005.}\\Department of Mathematics, Capital Normal University,\\zhaoxve@mail.cnu.edu.cn}
\date{}
\maketitle

\begin{abstract}
Let $G$ be a compact and $1$--connected Lie group with a maximal torus $T$.
Based on Schubert calculus on the flag manifold $G/T$ \cite{DZ4} we construct
the integral cohomology ring $H^{\ast}(G)$ uniformly for all $G$.

\begin{description}
\item[2000 Mathematical Subject Classification: ] 14M15; 55T10

\item[Key words:] Lie groups, Cohomology, Schubert calculus

\end{description}
\end{abstract}

\section{Introduction}

In this paper, the Lie groups $G$ under consideration are compact and
$1$--connected; the coefficient ring $\mathcal{R}$ for the cohomologies is
either the ring $\mathbb{Z}$ of integers, the field $\mathbb{R}$ of reals, or
one of the finite fields $\mathbb{F}_{p}$. For simplicity the symbol
$\mathbb{F}$ is also used to denote the field coefficients $\mathbb{R}$ or
$\mathbb{F}_{p}$. Write $H^{\ast}(X)$ instead of $H^{\ast}(X;\mathbb{Z})$ for
the integral cohomology of a topological space $X$.

The problem of computing the cohomology of Lie groups began with E. Cartan in
1929. It is a focus of algebraic topology for the fundamental role of Lie
groups playing in geometry and topology \cite[Chapter VI]{D},\cite{MT,Sa}.
Concerning the achievements by many mathematicians in about one century (e.g.
\cite{Br,Po,H,Y,L,BC,B1,B2,B3,AS,A2,P}) the following problem remains.

\bigskip

\noindent\textbf{Problem 1.1.} \textsl{Find a unified construction of the ring
}$H^{\ast}(G;\mathcal{R})$\textsl{ that is free of the types of the Lie group
}$G$\textsl{.}

\textsl{In particular, determine the integral cohomology ring }$H^{\ast}%
(G)$\textsl{ of the five exceptional Lie groups }$G=G_{2},F_{4},E_{6}%
,E_{7},E_{8}$\textsl{.}

\bigskip

\noindent In this paper we solve the problem by an explicit construction of
the cohomology $H^{\ast}(G;\mathcal{R})$ in the context of Schubert calculus
\cite{DZ4} for the $1$--connected Lie groups $G$. In the sequel work \cite{D2}
this construction will be extended to all compact Lie groups.

For a $1$--connected Lie group $G$ with a maximal torus $T$ let $\{\omega
_{1},\cdots,\omega_{n}\}$ $\subset H^{2}(G/T)$ be a set of fundamental
dominant weights, where $n=\dim T$ \cite{BH}. Our construction factors through
the Leray--Serre spectral sequence $\{E_{r}^{\ast,\ast}(G;\mathcal{R}%
),d_{r}\}$ of the fibration

\begin{enumerate}
\item[(1.1)] $T\hookrightarrow G\overset{\pi}{\rightarrow}G/T$\ \ 
\end{enumerate}

\noindent for which one has (\cite{M1}):

\begin{enumerate}
\item[(1.2)] $E_{2}^{\ast,\ast}(G;\mathcal{R})=H^{\ast}(G/T)\otimes
\Lambda_{\mathcal{R}}^{\ast}(t_{1},\cdots,t_{n})$,

\item[(1.3)] the differential $d_{2}:E_{2}^{p,q}(G;\mathcal{R})\rightarrow
E_{2}^{p+2,q-1}(G;\mathcal{R})$ is given by

$\quad d_{2}(x\otimes t_{k})=x\omega_{k}\otimes1$, $x\in H^{p}(G/T)$, $1\leq
k\leq n$.
\end{enumerate}

\noindent where $t_{i}\in H^{1}(T)$ is the class that is mapped to $\omega
_{i}$ under the transgression $\tau:H^{1}(T)\rightarrow H^{2}(G/T)$ \cite{BH},
and where $H^{\ast}(T;\mathcal{R})=\Lambda_{\mathcal{R}}^{\ast}(t_{1}%
,\cdots,t_{n})$ is the exterior ring generated by the $t_{i}$'s over
$\mathcal{R}$. To have the spectral sequence calculably we need a concise
characterization of the factor ring $H^{\ast}(G/T)$ in (1.2). The following
result from Schubert calculus serves this purpose.

\bigskip

\noindent\textbf{Theorem 1.2 (\cite[Theorem 1.2]{DZ4}). }\textsl{There exists
a set }$\left\{  y_{1},\cdots,y_{m}\right\}  $\textsl{\ of Schubert classes on
}$G/T$\textsl{\ with }$\deg y_{i}$\textsl{\ }$>2$,\textsl{\ so that the set
}$\{\omega_{1},\cdots,\omega_{n},y_{1},\cdots,$ $y_{m}\}$\textsl{\ is a
minimal system of generators of the ring }$H^{\ast}(G/T)$\textsl{.}

\textsl{Moreover, with respect to these generators one} \textsl{has the
presentation}

\begin{enumerate}
\item[(1.4)] $H^{\ast}(G/T)=\mathbb{Z}[\omega_{1},\cdots,\omega_{n}%
,y_{1},\cdots,y_{m}]/\left\langle h_{i},f_{j},g_{j}\right\rangle _{1\leq i\leq
k;1\leq j\leq m}$\textsl{,}
\end{enumerate}

\noindent\textsl{such that}

\begin{quote}
\textsl{i) for each }$1\leq i\leq k$\textsl{, }$h_{i}\in\left\langle
\omega_{1},\cdots,\omega_{n}\right\rangle $\textsl{;}

\textsl{ii) for each }$1\leq j\leq m$\textsl{, the pair }$(f_{j},g_{j}%
)$\textsl{\ of polynomials is related to the Schubert class }$y_{j}%
$\textsl{\ in the fashion}

$\qquad f_{j}$\textsl{\ }$=$\textsl{\ }$p_{j}y_{j}+\alpha_{j}$\textsl{, \quad
}$g_{j}=y_{j}^{k_{j}}+\beta_{j}$\textsl{,}

\noindent\textsl{with }$p_{j}\in\{2,3,5\}$\textsl{\ and }$\alpha_{j},\beta
_{j}\in\left\langle \omega_{1},\cdots,\omega_{n}\right\rangle $\textsl{.}%
$\square$
\end{quote}

Two remarks on Theorem 1.2 are in order. Firstly, with the minimum constraint
on the integers $k$ and $m$ the sets $\{\deg h_{i}\mid1\leq i\leq k\}$,
$\{\deg y_{j},p_{j},k_{j}\mid1\leq j\leq m\}$ of integers appearing in Theorem
1.2 can be shown to be invariants of the group $G$, and will be called
\textsl{the basic data }of $G$. For all the $1$--connected simple Lie groups
$G=SU(n+1),Sp(n),Spin(n+2),G_{2},F_{4},E_{6},E_{7}$ and $E_{8}$ their basic
data are tabulated in Section 5.

Next, for a simple Lie group $G\neq E_{8}$ the polynomials $h_{i},f_{j},g_{j}$
in (1.4) can be shown to be algebraically independent. In contrast, in term of
the basic data for the group $E_{8}$ given in Table 5.2 there appears the
following phenomena which will cause a few additional concern for the group
$E_{8}$ in our unified approach to the ring $H^{\ast}(G;\mathcal{R})$.

\bigskip

\noindent\textbf{Lemma 1.3 (\cite[Theorem 1.3]{DZ4}).} \textsl{For }$G=E_{8}%
$\textsl{ there exists a polynomial }$\phi$\textsl{ of the form }$\phi
=2y_{4}^{5}-y_{6}^{3}+y_{7}^{2}+\beta$\textsl{ with }$\beta\in\left\langle
\omega_{1},\cdots,\omega_{8}\right\rangle $\textsl{ so that}

\begin{enumerate}
\item[(1.5)] $\left\{
\begin{tabular}
[c]{l}%
$g_{4}=-12\phi+5y_{4}^{4}f_{4}-4y_{6}^{2}f_{6}+6y_{7}f_{7}$;\\
$g_{6}=-10\phi+4y_{4}^{4}f_{4}-3y_{6}^{2}f_{6}+5y_{7}f_{7}$;\\
$g_{7}=15\phi-6y_{4}^{4}f_{4}+5y_{6}^{2}f_{6}-7y_{7}f_{7}$.$\square$%
\end{tabular}
\ \ \ \ \right.  $
\end{enumerate}

Inputting the presentation (1.4) of the ring $H^{\ast}(G/T)$ into the formula
(1.2) of $E_{2}^{\ast,\ast}(G;\mathcal{R})$ and utilizing the polynomials
$h_{i},\alpha_{j},\beta_{j}$, one can construct a set of explicit generators
for the ring $H^{\ast}(G;\mathcal{R})$. To explain how this construction
proceeds we take the case $\mathcal{R}=\mathbb{Z}$ as an illustrative example.

The construction requires two ingredients. Firstly, each polynomial
$P\in\left\langle \omega_{1},\cdots,\omega_{n}\right\rangle $ admits a
decomposition of the form $P=p_{1}\omega_{1}+\cdots+p_{n}\omega_{n}$ which
gives rise to a map

\begin{enumerate}
\item[(1.6)] $\varphi:\left\langle \omega_{1},\cdots,\omega_{n}\right\rangle
\rightarrow E_{2}^{\ast,1}(G)=H^{\ast}(G/T)\otimes\Lambda_{\mathbb{Z}}%
^{1}(t_{1},\cdots,t_{n})$
\end{enumerate}

\noindent by $\varphi(P)=p_{1}\otimes t_{1}+\cdots+p_{n}\otimes t_{n}$ (see
Lemma 2.2). Secondly, since $E_{2}^{p,q}(G)=0$ for odd $p$ one has
the\textsl{\ canonical }homomorphism

\begin{enumerate}
\item[(1.7)] $\kappa:$ $E_{3}^{2k,1}(G)\rightarrow H^{2k+1}(G)$ (see (2.21))
\end{enumerate}

\noindent that interprets directly elements of $E_{3}^{2k,1}(G)$ as cohomology
classes of $G$.

In view of the fibration (1.1) we get from the Schubert classes $y_{i}$ on
$G/T$ the integral cohomology classes of the group $G$

\begin{quote}
$x_{s}:=\pi^{\ast}(y_{i})\in H^{\ast}(G)$, $s=\deg y_{i}$, $1\leq i\leq m$.
\end{quote}

Granted with the maps $\varphi$ and $\kappa$ in (1.6) and (1.7) the
polynomials $h_{i}$, $\alpha_{j}$, $\beta_{j}$ in (1.4) gives rise to the
following integral cohomology classes

\begin{quote}
$\varrho_{s}:=\kappa\lbrack\varphi(h_{i})]\in H^{\ast}(G),$ $s=\deg h_{i}-1$,
$1\leq i\leq k$;

$\varrho_{t}:=\kappa\lbrack\varphi(p_{j}\beta_{j}-y_{j}^{k_{j}}\alpha_{j})]\in
H^{\ast}(G),$ $t=\deg\beta_{j}-1$,
\end{quote}

\noindent where $[\gamma]\in E_{3}^{\ast,\ast}(G)$ denotes the cohomology
class of a $d_{2}$--cocycle $\gamma\in E_{2}^{\ast,\ast}(G)$, and where $1\leq
j\leq m$ with $j\neq4,7$ if $G=E_{8}$.

Finally, for a prime $p\in\{2,3,5\}$ and a multi-index $J\subset
\{1,\cdots,m\}$ with $p_{t}=p$, $t\in J$, we set $I=\{\deg\alpha_{t}\mid t\in
J\}$ and let

\begin{quote}
$\mathcal{C}_{I}:=\beta_{p}(%
{\displaystyle\prod\limits_{s\in I}}
\xi_{s-1})\in H^{\ast}(G)$,
\end{quote}

\noindent where $\xi_{s-1}=$ $\kappa^{\prime}[\varphi^{\prime}(\alpha_{t}%
)]\in$ $H^{\ast}(G;\mathbb{F}_{p})$ with $s=\deg\alpha_{t}$, $\kappa^{\prime}$
and $\varphi^{\prime}$ are the $\mathbb{F}_{p}$--analogue of the maps $\kappa$
and $\varphi$ in (1.6) and (1.7), respectively, and where $\beta_{p}%
:H^{r}(G;\mathbb{F}_{p})\rightarrow H^{r+1}(G)$ is the Bockstein homomorphism.

Using these three types $x_{s}$, $\varrho_{t}$ and $\mathcal{C}_{I}$ of
cohomology classes just described a unified presentation of the additive
cohomology $H^{\ast}(G)$ is given by Theorem 3.6. It implies that

i) the elements $\varrho_{t}$ have infinite order whose square free products
form a basis for the free part of the cohomology $H^{\ast}(G)$;

ii) the $p$--\textsl{primary component }$\sigma_{p}(G)$ of the ring $H^{\ast
}(G)$ is generated by the classes $x_{s}=\pi^{\ast}(y_{i})$ and $\mathcal{C}%
_{I}$ with $s=\deg y_{i}$ and $p_{i}=p$ for $i\in I$.

To determine the multiplicative relations among these explicitly constructed
generators we shall make crucial use of the spectral sequence, and a coherent
calculation in cohomologies with different coefficients. In particular, the
following theorem, that presents the cohomology rings of the five exceptional
Lie groups by these generators, will be established in Section 4, where the
notation used can be found at the beginning of Section 2.

\bigskip

\noindent\textbf{Theorem 1.4.} \textsl{For the five exceptional Lie groups we
have}

\begin{enumerate}
\item[(1.8)] $H^{\ast}(G_{2})=\Delta(\varrho_{3})\otimes\Lambda(\varrho
_{11})\oplus\sigma_{2}(G_{2})$\textsl{ with}

$\qquad\sigma_{2}(G_{2})=\mathbb{F}_{2}[x_{6}]^{+}/\left\langle x_{6}%
^{2}\right\rangle \otimes\Delta_{\mathbb{F}_{2}}(\varrho_{3})$\textsl{, }

\textsl{where the generators are subject to the following relations}

$\qquad\varrho_{3}^{2}=x_{6}$\textsl{,} $x_{6}\varrho_{11}=0$\textsl{.}

\item[(1.9)] $H^{\ast}(F_{4})=\Delta(\varrho_{3})\otimes\Lambda(\varrho
_{11},\varrho_{15},\varrho_{23})\oplus\sigma_{2}(F_{4})\oplus\sigma_{3}%
(F_{4})$\textsl{ with}

$\qquad\sigma_{2}(F_{4})=\mathbb{F}_{2}[x_{6}]^{+}/\left\langle x_{6}%
^{2}\right\rangle \otimes\Delta_{\mathbb{F}_{2}}(\varrho_{3})\otimes
\Lambda_{\mathbb{F}_{2}}(\varrho_{15},\varrho_{23})$\textsl{,}

$\qquad\sigma_{3}(F_{4})=\mathbb{F}_{3}[x_{8}]^{+}/\left\langle x_{8}%
^{3}\right\rangle \otimes\Lambda_{\mathbb{F}_{3}}(\varrho_{3},\varrho
_{11},\varrho_{15})$\textsl{,}

\textsl{where the generators are subject to the following relations}

$\qquad\varrho_{3}^{2}=x_{6}$\textsl{, } $x_{6}\varrho_{11}=0$\textsl{,
}$x_{8}\varrho_{23}=0$\textsl{.}

\item[(1.10)] $H^{\ast}(E_{6})=\Delta(\varrho_{3})\otimes\Lambda(\varrho
_{9},\varrho_{11},\varrho_{15},\varrho_{17},\varrho_{23})\oplus\sigma
_{2}(E_{6})\oplus\sigma_{3}(E_{6})$\textsl{ with}

$\qquad\sigma_{2}(E_{6})=\mathbb{F}_{2}[x_{6}]^{+}/\left\langle x_{6}%
^{2}\right\rangle \otimes\Delta_{\mathbb{F}_{2}}(\varrho_{3})\otimes
\Lambda_{\mathbb{F}_{2}}(\varrho_{9},\varrho_{15},\varrho_{17},\varrho_{23}%
)$\textsl{,}

$\qquad\sigma_{3}(E_{6})=\mathbb{F}_{3}[x_{8}]^{+}/\left\langle x_{8}%
^{3}\right\rangle \otimes\Lambda_{\mathbb{F}_{3}}(\varrho_{3},\varrho
_{9},\varrho_{11},\varrho_{15},\varrho_{17})$\textsl{,}

\textsl{where the generators are subject to the following relations}

$\qquad\varrho_{3}^{2}=x_{6}$\textsl{, }$x_{6}\varrho_{11}=0$\textsl{, }%
$x_{8}\varrho_{23}=0$\textsl{.}

\item[(1.11)] $H^{\ast}(E_{7})=\Delta(\varrho_{3})\otimes\Lambda(\varrho
_{11},\varrho_{15},\varrho_{19},\varrho_{23},\varrho_{27},\varrho
_{35})\underset{p=2,3}{\oplus}\sigma_{p}(E_{7})$\textsl{ with}

$\qquad\sigma_{2}(E_{7})=\frac{\mathbb{F}_{2}[x_{6},x_{10},x_{18}%
,\mathcal{C}_{I}]^{+}}{\left\langle x_{6}^{2},x_{10}^{2},x_{18}^{2}%
,\mathcal{D}_{I},\mathcal{R}_{J},\mathcal{S}_{K,L},\mathcal{H}_{r,I}%
\right\rangle }\otimes\Delta_{\mathbb{F}_{2}}(\varrho_{3})\otimes
\Lambda_{\mathbb{F}_{2}}(\varrho_{15},\varrho_{23},\varrho_{27})$

\textsl{for }$I,J,K,L\subseteq\{6,10,18\}$\textsl{, }$\left\vert I\right\vert
,\left\vert J\right\vert ,\left\vert K\right\vert \geq2$\textsl{, }%
$r\in\{11,19,35\}$\textsl{;}

$\qquad\sigma_{3}(E_{7})=\frac{\mathbb{F}_{3}[x_{8}]^{+}}{\left\langle
x_{8}^{3}\right\rangle }\otimes\Lambda_{\mathbb{F}_{3}}(\varrho_{3}%
,\varrho_{11},\varrho_{15},\varrho_{19},\varrho_{27},\varrho_{35})$\textsl{,}

\textsl{where the generators are subject to the following relations}

$\qquad\varrho_{3}^{2}=x_{6}$\textsl{, }$x_{8}\varrho_{23}=0.$

\item[(1.12)] $H^{\ast}(E_{8})=\Delta(\varrho_{3},\varrho_{15},\varrho
_{23})\otimes\Lambda(\varrho_{27},\varrho_{35},\varrho_{39},\varrho
_{47},\varrho_{59})\underset{p=2,3,5}{\oplus}\sigma_{p}(E_{8})$\textsl{ with}

$\qquad\sigma_{2}(E_{8})=\frac{\mathbb{F}_{2}[x_{6},x_{10},x_{18}%
,x_{30},\mathcal{C}_{I}]^{+}}{\left\langle x_{6}^{8},x_{10}^{4},x_{18}%
^{2},x_{30}^{2},\mathcal{D}_{I},\mathcal{R}_{J},\mathcal{S}_{K,L}%
,\mathcal{H}_{r,I}\right\rangle }\otimes\Delta_{\mathbb{F}_{2}}(\varrho
_{3},\varrho_{15},\varrho_{23})\otimes\Lambda_{\mathbb{F}_{2}}(\varrho_{27})$

\textsl{for} $I$\textsl{,}$J,K\subseteq\{6,10,18,30\}$\textsl{, }$\left\vert
I\right\vert ,\left\vert J\right\vert ,\left\vert K\right\vert \geq2$\textsl{,
}$r\in\{35,39,47,59\}$\textsl{;}

$\qquad\sigma_{3}(E_{8})=\frac{\mathbb{F}_{3}[x_{8},x_{20},\mathcal{C}%
_{\{8,20\}}]^{+}}{\left\langle x_{8}^{3},x_{20}^{3},x_{8}^{2}x_{20}%
^{2}\mathcal{C}_{\{8,20\}},(\mathcal{C}_{\{8,20\}})^{2}\right\rangle }%
\otimes\Lambda_{\mathbb{F}_{3}}(\varrho_{3},\varrho_{15},\varrho_{27}%
,\varrho_{35},\varrho_{39},\varrho_{47})$\textsl{,}

$\qquad\sigma_{5}(E_{8})=\frac{\mathbb{F}_{5}[x_{12}]^{+}}{\left\langle
x_{12}^{5}\right\rangle }\otimes\Lambda_{\mathbb{F}_{5}}(\varrho_{3}%
,\varrho_{15},\varrho_{23},\varrho_{27},\varrho_{35},\varrho_{39},\varrho
_{47})$\textsl{,}

\textsl{where the generators are subject to the following relations}

$\qquad\varrho_{3}^{2}=x_{6}$\textsl{, }$\varrho_{15}^{2}=x_{30}$\textsl{,
}$\varrho_{23}^{2}=x_{6}^{6}x_{10}$\textsl{, }$x_{2s}\varrho_{3s-1}=0$
\textsl{for} $s=4,5$\textsl{,}

$\qquad x_{8}\varrho_{59}=x_{20}^{2}\mathcal{C}_{\{8,20\}}$\textsl{, }%
$x_{20}\varrho_{23}=x_{8}^{2}\mathcal{C}_{\{8,20\}}$\textsl{, }$x_{12}%
\varrho_{59}=0$\textsl{.}
\end{enumerate}

\noindent\textsl{In addition, the relations of the types} $\mathcal{D}%
_{J},\mathcal{R}_{K},\mathcal{S}_{I,J}$ \textsl{and} $\mathcal{H}_{r,I}%
$\textsl{,} \textsl{that occur only} \textsl{in the formulae of the ideals
}$\sigma_{2}(E_{7})$ \textsl{and} $\sigma_{2}(E_{8})$\textsl{, are given by
the formulae (4.4), (4.7) and (4.9), respectively.}$\square$

\bigskip

One may compare the presentations in (1.8) and (1.9) with the computation of
the rings $H^{\ast}(G_{2})$ and $H^{\ast}(F_{4})$ by Borel \cite{B2,B3}.

The remaining sections are so arranged. Section 2 develops preliminary results
on Koszul complexes, which are applied in Section 3 to describe the cohomology
$H^{\ast}(G;\mathcal{R})$ as an additive group. The calculation is extended in
Sections 4 as to determine the ring structure on $H^{\ast}(G;\mathcal{R})$ for
the exceptional Lie groups $G$.

Historically, Schubert calculus is the intersection theory of the 19th
century, together with its applications to enumerative geometry \cite{Sch}.
Clarifying this calculus was the content of Hilbert's problem 15
\cite{F,Kl,So}. van der Waerden and A. Weil, who secured the foundation of
modern intersection theory, attributed this calculus to the determination of
the cohomology rings of the flag manifolds \cite{Wa}, \cite[p.331]{W}.

The idea to study the cohomology of Lie groups by applying Schubert calculus
on the flag manifolds $G/T$ began with Ka\v{c} \cite{K1,K2} and Marlin
\cite{M1}. They conjectured respectively that the integral Leray--Serre
spectral sequence of the fibration (1.1) satisfies the following relations

\begin{quote}
$E_{3}^{\ast,\ast}(G)=E_{\infty}^{\ast,\ast}(G)$, $E_{3}^{\ast,\ast
}(G)=H^{\ast}(G)$.
\end{quote}

\noindent These conjectures will be confirmed in Section 3, see Remark 3.7.

\section{Preliminaries in Koszul complexes}

We fix certain notation concerning graded modules or rings. For a graded ring
$A$ and a finite set $S=\{u_{1},\cdots,u_{t}\}$ of graded elements we write

\begin{quote}
$A\{S\}=A\{u_{i}\}_{1\leq i\leq t}$ for the free $A$--module with basis
$\{u_{1},\cdots,u_{t}\}$;

$A[S]=A[u_{i}]_{1\leq i\leq t}$ for the polynomials ring in the $u_{i}$'s over
$A$;

$\Lambda_{\mathcal{R}}(S)=\Lambda_{\mathcal{R}}(u_{i})_{1\leq i\leq t}$ for
the exterior ring in the $u_{i}$'s over $\mathcal{R}$;

$A\otimes\Delta(S)=A\otimes\Delta(u_{i})_{1\leq i\leq t}$ for the free
$A$--module in the simple system $\left\{  u_{1},\cdots,u_{t}\right\}  $ of
generators \cite{B4}.
\end{quote}

\noindent If $A=\mathbb{Z}$ (resp. $\mathbb{F}_{p}$) we use $\Delta(S)$ (resp.
$\Delta_{\mathbb{F}_{p}}(S)$) instead of $A\otimes\Delta(S)$.

For a finite subset $\{a_{1},\cdots,a_{m}\}$ of a ring $A$ let $\left\langle
a_{1},\cdots,a_{m}\right\rangle $ be the ideal generated by $a_{1}%
,\cdots,a_{m}$, and write $A/\left\langle a_{1},\cdots,a_{m}\right\rangle $
for the quotient ring.

If $V=V^{0}\oplus V^{1}\oplus\cdots$ is a graded vector space (resp. ring),
define its subspace (resp. subring) $V^{+}$ by $V^{1}\oplus V^{2}\oplus\cdots
$. For instance if $\mathbb{Z}[x_{1},\cdots,x_{n}]$ is the ring of integral
polynomials with $\deg x_{i}>0$, then $\mathbb{Z}[x_{1},\cdots,x_{n}]^{+}$ is
the subring of the polynomials without constant terms.

All elements $a,b,\cdots$ in a graded module or ring will assumed to be
homogeneous, whose degrees are denoted by $\left\vert a\right\vert ,\left\vert
b\right\vert ,\cdots$.

Given a ring $A$ and an ordered subset $\{x_{1},\cdots,x_{n}\}$ \textsl{the
Koszul complex} of $A$ with respect to $\left\{  x_{1},\cdots,x_{n}\right\}  $
is the complex $\{C(A);\delta\}$ defined by

\begin{quote}
i) $C(A)=A\otimes\Delta(\theta_{1},\cdots,\theta_{n})$;

ii) $\delta(1\otimes\theta_{i})=x_{i}\otimes1$; $\delta(a\otimes1)=0$, $a\in
A$.
\end{quote}

\noindent Let $\Delta^{k}$ be the subgroup of $\Delta(\theta_{1},\cdots
,\theta_{m})$ spanned by the square free monomials in $\theta_{1}%
,\cdots,\theta_{m}$ with length $k$. Then the $k^{th}$ cohomology group of the
complex $\{C(A);\delta\}$ is

\begin{quote}
$H^{k}(\{C(A);\delta\})=\frac{\ker[\delta:A\otimes\Delta^{k}\rightarrow
A\otimes\Delta^{k-1}]}{\operatorname{Im}[\delta:A\otimes\Delta^{k+1}%
\rightarrow A\otimes\Delta^{k}]}$.
\end{quote}

\noindent In particular, the $E_{2}$--page of the spectral sequence of the
fibration (1.1) is \textsl{the Koszul complex} of the ring $H^{\ast
}(G/T;\mathcal{R})$ with respect to the set $\{\omega_{1},\cdots,\omega_{n}\}$
of weights whose cohomology is the third page $E_{3}^{\ast,\ast}%
(G;\mathcal{R})$. In this section we develop some general results on the
cohomology of certain Koszul complexes that are useful for us to formulate
$E_{3}^{\ast,\ast}(G;\mathcal{R})$ (resp. $H^{\ast}(G;\mathcal{R})$) in term
of the presentation (1.4) on the ring $H^{\ast}(G/T)$.

\subsection{The cohomology of a Koszul complex}

\noindent Given a graded ring $A$ with the following presentation

\begin{enumerate}
\item[(2.1)] $A=\mathbb{F}_{p}[y_{1},\cdots,y_{r}]/\left\langle y_{1}^{k_{1}%
},\cdots,y_{r}^{k_{r}}\right\rangle $, $\left\vert y_{t}\right\vert
\equiv0\operatorname{mod}2$, $1\leq t\leq r$,
\end{enumerate}

\noindent let $\left\{  C(A),\delta\right\}  $ be the \textsl{Koszul complex}
of $A$ with respect to the ordered subset $(-y_{1},\cdots,-y_{r})$. In view of
the formula $C(A)=A\otimes\Delta(\theta_{1},\cdots,\theta_{r})$ we introduce
for each multi--index $I\subseteq\{1,\cdots,r\}$ the following elements of
$C(A)$

\begin{quote}
$\theta_{I}=\prod\limits_{t\in I}\theta_{t}$; $\ f_{I}=\sum\limits_{t\in
I}-y_{t}\theta_{I_{t}}(=\delta(\theta_{I}))$; $\quad g_{I}=(\prod\limits_{t\in
I}y_{t}^{k_{t}-1})\theta_{I}$,
\end{quote}

\noindent where $I_{t}$ denotes the complement of $t$ in $I$. We shall also set

\begin{enumerate}
\item[(2.2)] $c_{I}=\delta\theta_{I}$, $\quad D_{I}=\sum\limits_{t\in I}%
-y_{t}c_{I_{t}}$, $\quad R_{I}=(\prod\limits_{t\in I}y_{t}^{k_{t}-1})c_{I}$.
\end{enumerate}

\noindent\textbf{Theorem 2.1.} \textsl{For the Koszul complex }$\left\{
C(A),\delta\right\}  $ \textsl{we have}

\begin{quote}
\textsl{i)} $H^{\ast}(\left\{  C(A),\delta\right\}  )=\Delta_{\mathbb{F}_{p}%
}(y_{1}^{k_{1}-1}\theta_{1},\cdots,y_{r}^{k_{r}-1}\theta_{r})$\textsl{;}

\textsl{ii)} $\operatorname{Im}\delta=\frac{A\{1,c_{I}\}^{+}}{\left\langle
D_{J},R_{K}\right\rangle }$,\textsl{\ }$I,J,K\subseteq\{1,\cdots
,r\}$\textsl{\ with }$\left\vert I\right\vert ,\left\vert J\right\vert
,\left\vert K\right\vert \geq2$\textsl{,\ }
\end{quote}

\noindent\textsl{where }$\left\langle D_{J},R_{K}\right\rangle \subset
A\{1,c_{I}\}^{+}$ \textsl{denotes the} \textsl{sub} $A$\textsl{--module}
\textsl{spanned by} \textsl{the elements} $D_{J}$ \textsl{and} $R_{K}%
$\textsl{.}

\bigskip

\noindent\textbf{Proof.} Since $C(A)=$ $\underset{1\leq t\leq r}{\otimes}%
C_{t}$ with $C_{t}=(\mathbb{F}_{p}[y_{t}]/\left\langle y_{t}^{k_{t}%
}\right\rangle )\otimes\Lambda_{\mathbb{F}_{p}}(\theta_{t})$ a closed subspace
of $\delta$, one has by the K\"{u}nneth formula that

\begin{quote}
$H^{\ast}(\left\{  C(A),\delta\right\}  )=\underset{1\leq t\leq r}{\otimes
}H^{\ast}(\mathcal{C}_{t},\delta)$
\end{quote}

\noindent Formula i) is shown by the obvious fact that the cohomology
$H^{\ast}(C_{t},\delta)$ has a basis represented by the subset $\left\{
1,y_{t}^{k_{t}-1}\theta_{t}\right\}  \subset C_{t}$, $1\leq t\leq r$.

In view of the short exact sequence $0\rightarrow\ker\delta\rightarrow
C(A)\overset{\delta}{\rightarrow}\operatorname{Im}\delta\rightarrow0$ the
differential $\delta$ induces an isomorphism of $A$--modules

\begin{quote}
a) $\overline{\delta}:C(A)/\ker\delta\overset{\cong}{\rightarrow
}\operatorname{Im}\delta$.
\end{quote}

\noindent In the quotient space $C(A)/\ker\delta$ the numerator $C(A)$ has the
$A$ basis $\{1,\theta_{I}\mid I\subseteq\{1,\cdots,r\}\}$, while the
denominator admits the decomposition

\begin{quote}
$\ker\delta=H^{\ast}(\left\{  C(A),\delta\right\}  )\oplus\operatorname{Im}%
\delta$
\end{quote}

\noindent in which the first summand $H^{\ast}(\left\{  C(A),\delta\right\}
)$ has the $\mathbb{F}_{p}$--basis $\{1,g_{K}\mid K\subseteq\{1,\cdots,r\}\}$
by i), and the second summand $\operatorname{Im}\delta$ is spanned over $A$ by
the subset $\{f_{J}=\delta(\theta_{J})\mid J\subseteq\{1,\cdots,r\}\}$. That is

\begin{quote}
b) $C(A)/\ker\delta=\frac{A\{1,\theta_{I}\}^{+}}{\left\langle f_{J}%
,g_{K}\right\rangle },I,J,K\subseteq\{1,\cdots,r\}$.
\end{quote}

\noindent The presentation ii) is shown by a) and b), the obvious relations

\begin{quote}
$\overline{\delta}(\theta_{I})=c_{I}$, $\overline{\delta}(f_{J})=D_{J}$,
$\overline{\delta}(g_{K})=R_{K}$,
\end{quote}

\noindent together with the facts that if $I=\{t\}$ is a singleton, then
$c_{I}=-y_{t}$, $D_{I}=0$, $R_{I}=-y_{t}^{k_{t}}=0$.$\square$

\subsection{The Koszul complex of a regularly truncated polynomial ring}

Let $B$ be a graded ring with the following presentation

\begin{enumerate}
\item[(2.3)] $B=\frac{\mathcal{R}[\omega_{1},\cdots,\omega_{n},y_{1}%
,\cdots,y_{t}]}{\left\langle r_{1},\cdots,r_{m}\right\rangle }$, $2=\left\vert
\omega_{i}\right\vert <\left\vert y_{1}\right\vert <\cdots<\left\vert
y_{t}\right\vert $,
\end{enumerate}

\noindent and let $\left\{  C(B),\delta\right\}  $ be the Koszul complex of
$B$ with respect to the ordered subset $\left\{  \omega_{1},\cdots,\omega
_{n}\right\}  $ with

\begin{quote}
i) $C(B)=B\otimes\Lambda(t_{1},\cdots,t_{n})$, $\left\vert t_{i}\right\vert
=1$;

ii) $\delta(1\otimes t_{i})=\omega_{i}\otimes1$, $\delta(a\otimes1)=0$, $1\leq
i\leq n,a\in B$.
\end{quote}

\noindent For a polynomial $P\in\left\langle \omega_{1},\cdots,\omega
_{n}\right\rangle \subset\mathcal{R}[\omega_{1},\cdots,\omega_{n},y_{1}%
,\cdots,y_{t}]$ we can write

\begin{enumerate}
\item[(2.4)] $P=p_{1}\cdot\omega_{1}+\cdots+p_{n}\cdot\omega_{n}$ with
$p_{i}\in\mathcal{R}[\omega_{i},y_{j}]$.
\end{enumerate}

\noindent and set $\widetilde{P}:=p_{1}\otimes t_{1}+\cdots+$ $p_{n}\otimes
t_{n}\in B\otimes\Lambda^{1}$. It is crucial to notice that if $P\in
\left\langle \omega_{1},\cdots,\omega_{n}\right\rangle \cap\left\langle
r_{1},\cdots,r_{m}\right\rangle $ then $\widetilde{P}\in\ker\delta$. For a
$\delta$--cocycle $\gamma\in C(B)$ let $\left[  \gamma\right]  \in H^{\ast
}(\left\{  C(B),\delta\right\}  )$ be its cohomology class.

\bigskip

\noindent\textbf{Lemma 2.2.} \textsl{The map }$\varphi:\left\langle \omega
_{1},\cdots,\omega_{n}\right\rangle \cap\left\langle r_{1},\cdots
,r_{m}\right\rangle \rightarrow H^{1}(\left\{  C(B),\delta\right\}
)$\textsl{\ by }$\varphi(P)=[\widetilde{P}]$\textsl{\ is well defined. }

\bigskip

\noindent\textbf{Proof.}\textsl{\ }We are bound to show that the class
$[\widetilde{P}]\in H^{1}(\left\{  C(B),\delta\right\}  )$ with $P\in
\left\langle \omega_{1},\cdots,\omega_{n}\right\rangle \cap\left\langle
r_{1},\cdots,r_{m}\right\rangle $ is independent of a choice of the expansion
(2.4). Assume in addition to (2.4) one has a second decomposition

\begin{quote}
$P=h_{1}\cdot\omega_{1}+\cdots+h_{n}\cdot\omega_{n}$.
\end{quote}

\noindent Then the equation

\begin{quote}
$(p_{1}-h_{1})\cdot\omega_{1}+\cdots+(p_{n}-h_{n})\cdot\omega_{n}=0$
\end{quote}

\noindent holds in the ring $\mathcal{R}[\omega_{i},y_{j}]$. We can assume
below that $p_{1}-h_{1}\neq0$.

Since the set $\left\{  \omega_{1},\cdots,\omega_{n}\right\}  $\ is
algebraically independent in the free polynomial ring $\mathcal{R}[\omega
_{i},y_{j}]$ the above equation implies that all the differences $p_{i}-h_{i}$
with $i\neq1$ are divisible by $\omega_{1}$. That is $p_{i}-h_{i}=q_{i}%
\cdot\omega_{1}$ for some $q_{i}\in\mathbb{Z}[\omega_{i},y_{j}]$, $2\leq i\leq
n$. The lemma is shown by the calculation

\begin{quote}
$\delta((q_{2}\otimes t_{1}t_{2}+\cdots+q_{n})\otimes t_{1}t_{n})=(p_{1}%
-h_{1})\otimes t_{1}+\cdots+(p_{n}-h_{n})\otimes t_{n}$.$\square$
\end{quote}

\bigskip

\noindent\textbf{Definition 2.3.} Let $B$ be a ring over a field
$\mathcal{R}=\mathbb{F}$ with the presentation (2.3). We call $B$ a
\textsl{regularly truncated} \textsl{polynomial ring} if the following two
additional conditions are satisfied

i) $n+t=m$;

ii) the set $\left\{  r_{1},\cdots,r_{n+t}\right\}  $ of relations admits a
partition $\left\{  \alpha_{1},\cdots,\alpha_{n}\right\}  \amalg\left\{
\gamma_{1},\cdots,\gamma_{t}\right\}  $ such that

\begin{quote}
a) $\alpha_{i}\in\left\langle \omega_{1},\cdots,\omega_{n}\right\rangle $,
$1\leq i\leq n$;

b) $\gamma_{j}=y_{j}^{k_{j}}+\beta_{j}$ with $\beta_{j}\in\left\langle
\omega_{1},\cdots,\omega_{n}\right\rangle $, $1\leq j\leq t$.$\square$
\end{quote}

Let $\left\{  C(B),\delta\right\}  $ be the Koszul complex associated to the
regularly truncated polynomial ring $B$. By the properties a) and b) one has that

\begin{quote}
$H^{0}(\left\{  C(B),\delta\right\}  )=B/\left\langle \omega_{1},\cdots
,\omega_{n}\right\rangle =\frac{\mathbb{F}[y_{1},\cdots,y_{t}]}{\left\langle
y_{1}^{k_{1}},\cdots,y_{t}^{k_{t}}\right\rangle }$.
\end{quote}

\noindent In addition, the construction $\varphi$ in Lemma 2.2 yields the
cohomology classes

\begin{quote}
$\xi_{\alpha_{i}}:=\varphi\left[  \widetilde{\alpha}_{i}\right]  \in
H^{1}(\left\{  C(B),\delta\right\}  )$, $1\leq i\leq n$.
\end{quote}

\noindent Since $\xi^{2}=0$ for any $\xi\in H^{1}(\left\{  C(B),\delta
\right\}  )$ (with respect to the product on $H^{\ast}(\left\{  C(B),\delta
\right\}  )$ inherited from that on $C(B)$) the inclusions $H^{r}(\left\{
C(B),\delta\right\}  )$ $\subset H^{\ast}(\left\{  C(B),\delta\right\}  )$
with $r=0,1$ extend naturally to a ring map

\begin{enumerate}
\item[(2.5)] $f:\frac{\mathbb{F}[y_{1},\cdots,y_{t}]}{\left\langle
y_{1}^{k_{1}},\cdots,y_{t}^{k_{t}}\right\rangle }\otimes\Lambda(\xi
_{\alpha_{i}})_{1\leq i\leq n}\rightarrow H^{\ast}(\left\{  C(B),\delta
\right\}  )$.
\end{enumerate}

\noindent\textbf{Theorem 2.4. }\textsl{If }$B$ \textsl{is a regularly
truncated} \textsl{polynomial ring with }$\dim_{\mathbb{F}}B<\infty$\textsl{,}
\textsl{then} \textsl{the map }$f$\textsl{ is surjective.}

\bigskip

\noindent\textbf{Proof of Theorem 2.4.} Set $B_{1}=\mathbb{F}[\omega
_{1},\cdots,\omega_{n},y_{1},\cdots,y_{t}]/\left\langle \gamma_{1}%
,\cdots,\gamma_{t}\right\rangle $ and let $r:B_{1}\rightarrow B$ be the
obvious quotient map. Consider the short exact sequence of Koszul complexes

\begin{enumerate}
\item[(2.6)] $0\rightarrow\left\{  C(\ker r),\delta_{2}\right\}
\overset{j^{\ast}}{\rightarrow}\left\{  C(B_{1}),\delta_{1}\right\}
\overset{r^{\ast}}{\rightarrow}\left\{  C(B),\delta\right\}  \rightarrow0$.
\end{enumerate}

\noindent associated to the short exact sequence $0\rightarrow\ker
r\overset{j}{\rightarrow}B_{1}\overset{r}{\rightarrow}B\rightarrow0$ of rings.
Since $B_{1}$ is a free $\mathbb{F}[\omega_{i}]$ module with the basis
$\left\{  y_{1}^{r_{1}}\cdots y_{t}^{r_{t}}\mid0\leq r_{i}<k_{i}\right\}  $ by
b), we have in addition to $H^{k}(\left\{  C(B_{1}),\delta_{1}\right\}  )=0$
($k\geq1$) that the map $r^{\ast}$ in (2.6) induces an isomorphism

\begin{quote}
$H^{0}(\left\{  C(B_{1}),\delta_{1}\right\}  )\cong H^{0}(\left\{
C(B),\delta\right\}  )$ ($=\mathbb{F}[y_{1},\cdots,y_{t}]/\left\langle
y_{1}^{k_{1}},\cdots,y_{t}^{k_{t}}\right\rangle $).
\end{quote}

\noindent In particular, the connecting homomorphism $\partial$ in the
cohomology exact sequence associated to (2.6) yields the following isomorphism
of $\mathbb{F}$--spaces, which is the key formula for us to establish the theorem:

\begin{enumerate}
\item[(2.7)] $\partial:H^{k}(\left\{  C(B),\delta\right\}  )\cong
H^{k-1}(\left\{  C(\ker r),\delta_{2}\right\}  )$, $k\geq1$.
\end{enumerate}

Since the number $n+t$ of the generators $\omega_{i},y_{j}$ equals to the
number $m$ of the relations by i) of Definition 2.3, the property
$\dim_{\mathbb{F}_{p}}B<\infty$ implies that the ordered sequence $\left\{
\gamma_{1},\cdots,\gamma_{t},\alpha_{1},\cdots,\alpha_{n}\right\}  $ is a
regular sequence of the ring $\mathbb{F}[\omega_{i},y_{j}]$ \cite[p.423]{E}.
As a result $\left\{  \alpha_{1},\cdots,\alpha_{n}\right\}  $ is a regular
sequence of $B_{1}$. By \cite[Corollary 17.5]{E} the ideal $\ker r$, as a
module over $B_{1}$, admits the following presentation

\begin{enumerate}
\item[(2.8)] $\ker r=$ $\frac{B_{1}\cdot\{x_{1},\cdots,x_{n}\}}{B_{1}%
\cdot\left\{  \alpha_{i}\cdot x_{j}-\alpha_{j}\cdot x_{i},1\leq i<j\leq
n\right\}  }$,
\end{enumerate}

\noindent where $\{x_{1},\cdots,x_{n}\}$ is the set of indeterminacies that is
characterized by the relations $j(b\cdot x_{i})=b\alpha_{i}$, $b\in B_{1}$,
$1\leq i\leq n$.

Let $\Lambda_{\alpha}^{k}\subset$ $\Lambda(\xi_{\alpha_{i}})_{1\leq i\leq n}$
(resp. $\Lambda^{k}\subset\Lambda(t_{i})_{1\leq i\leq n}$) be the subspace
spanned by the square free monomials in the $\xi_{\alpha_{i}}$'s (resp. in the
$t_{i}$'s) with length $k$, and let $f_{k}$ be the restriction of $f$ to the
subspace $\frac{\mathbb{F}[y_{1},\cdots,y_{t}]}{\left\langle y_{1}^{k_{1}%
},\cdots,y_{t}^{k_{t}}\right\rangle }\otimes\Lambda_{\alpha}^{k}$. It suffices
to show, by an induction on $k\geq1$, that $f_{k}$ is onto $H^{k}(\left\{
C(B),\delta\right\}  )$. If $k=1$ the result comes from (2.7), together with
the obvious relations

\begin{quote}
$\partial(\xi_{\alpha_{i}})=x_{i}$; $H^{0}(\left\{  C(\ker r),\delta
_{2}\right\}  )=\frac{\mathbb{F}[y_{1},\cdots,y_{t}]}{\left\langle
y_{1}^{k_{1}},\cdots,y_{t}^{k_{t}}\right\rangle }\cdot\{x_{1},\cdots,x_{n}\}$.
\end{quote}

For the remaining cases $k\geq2$ we examine the short exact ladder

\begin{quote}%
\begin{tabular}
[c]{lllllll}%
$0\rightarrow$ & $\ker r\otimes\Lambda^{k}$ & $\overset{j^{\ast}}{\rightarrow
}$ & $B_{1}\otimes\Lambda^{k}$ & $\overset{r^{\ast}}{\rightarrow}$ &
$B\otimes\Lambda^{k}$ & $\rightarrow0$\\
& $\delta_{2}\downarrow$ &  & $\delta_{1}\downarrow$ &  & $\delta\downarrow$ &
\\
$0\rightarrow$ & $\ker r\otimes\Lambda^{k-1}$ & $\overset{j^{\ast}%
}{\rightarrow}$ & $B_{1}\otimes\Lambda^{k-1}$ & $\overset{r^{\ast}%
}{\rightarrow}$ & $B\otimes\Lambda^{k-1}$ & $\rightarrow0$%
\end{tabular}

\end{quote}

\noindent induced by $\delta$. For a $z\in B_{1}\otimes\Lambda^{k}$ with
$r^{\ast}(z)\in\ker\delta$ one infers from

\begin{quote}
$\delta_{1}(z)=c_{1}\alpha_{1}+\cdots+c_{n}\alpha_{n}$, $c_{i}\in B_{1}%
\otimes\Lambda^{k-1}$
\end{quote}

\noindent that $j^{\ast-1}(\delta_{1}(z))=c_{1}\cdot x_{1}+\cdots+c_{n}\cdot
x_{n}\in\ker r\otimes\Lambda^{k-1}$. Consequently

\begin{enumerate}
\item[(2.9)] $\delta_{2}(j^{\ast-1}(\delta_{1}(z)))=\delta_{1}(c_{1})\cdot
x_{1}+\cdots+\delta_{1}(c_{n})\cdot x_{n}$.
\end{enumerate}

On the other hand since the element $j^{\ast-1}(\delta_{1}(z))$ is $\delta
_{2}$ closed, by (2.8) there exist unique $g_{i,j}\in B_{1}\otimes
\Lambda^{k-1}$ so that

\begin{enumerate}
\item[(2.10)] $\delta_{2}(j^{\ast-1}(\delta_{1}(z)))=%
{\displaystyle\sum\limits_{1\leq i<j\leq n}}
g_{i,j}(\alpha_{i}\cdot x_{j}-\alpha_{j}\cdot x_{i})$.
\end{enumerate}

\noindent Comparing this with (2.9) yields the equations

\begin{enumerate}
\item[(2.11)] $\delta_{1}(c_{i})=g_{1,i}\alpha_{1}+\cdots+g_{i-1,i}%
\alpha_{i-1}-g_{i,i+1}\alpha_{i+1}-\cdots-g_{i,n}\alpha_{n}$
\end{enumerate}

\noindent indicating that $r^{\ast}(c_{i})\in\ker\delta$. In addition, with
$\delta_{2}^{2}=0$ the formula (2.10) tells also that $\delta_{1}(g_{i,j})=0$.
By the inductive hypothesis we can assume that

\begin{enumerate}
\item[(2.12)] $\left[  r^{\ast}(g_{i,j})\right]  \in\operatorname{Im}f_{k-2}$,
$\left[  r^{\ast}(c_{i})\right]  \in\operatorname{Im}f_{k-1}$.
\end{enumerate}

Finally, consider the class

\begin{quote}
$z^{\prime}=%
{\displaystyle\sum\limits_{1\leq i\leq n}}
c_{i}\wedge\xi_{\alpha_{i}}-%
{\displaystyle\sum\limits_{1\leq i<j\leq n}}
g_{i,j}\wedge(\xi_{\alpha_{i}}\wedge\xi_{\alpha_{j}})\in B_{1}\otimes
\Lambda^{k}$.
\end{quote}

\noindent From $\delta_{1}(z^{\prime})=\delta_{1}(z)$ by (2.11) one gets that
$\partial\lbrack r^{\ast}(z)]=\partial\lbrack r^{\ast}(z^{\prime})]$. The
proof is completed by the following calculation due to (2.7) and (2.12)

\begin{quote}
$[r^{\ast}(z)]=[r^{\ast}(z^{\prime})]$

$=%
{\displaystyle\sum\limits_{1\leq i\leq n}}
\left[  r^{\ast}(c_{i})\right]  \wedge f_{1}(\xi_{\alpha_{i}})-%
{\displaystyle\sum\limits_{1\leq i<j\leq n}}
\left[  r^{\ast}(g_{i,j})\right]  \wedge f_{2}(\xi_{\alpha_{i}}\wedge
\xi_{\alpha_{j}})$.$\square$
\end{quote}

Theorem 2.4 is useful to formulate the term $E_{3}^{\ast,\ast}(G;\mathbb{F})$
out of the presentation (1.4) on the ring $H^{\ast}(G/T;\mathbb{F})$ because
of the following result.

\bigskip

\noindent\textbf{Lemma 2.5.} \textsl{For a simple Lie group }$G$\textsl{ with
a maximal torus }$T$\textsl{, the cohomology }$H^{\ast}(G/T;\mathbb{F}%
)$\textsl{ is a regularly truncated} \textsl{polynomial ring.}

\bigskip

\noindent\textbf{Proof. }Assume below that the rank of the simple Lie group
$G$ is $n$. Since the integral cohomology $H^{\ast}(G/T)$ is torsion free
\cite{BS}, a presentation of the cohomology $H^{\ast}(G/T;\mathbb{F})$ can be
deduced from (1.4), together with the isomorphism $H^{\ast}(G/T;\mathbb{F}%
)=H^{\ast}(G/T)\otimes\mathbb{F}$.

If $\mathbb{F}=\mathbb{R}$ we have $y_{j}=-\frac{1}{p_{j}}\alpha_{j}$ by the
relation of the type $f_{j}$ in Theorem 1.2. It implies that

\begin{enumerate}
\item[(2.13)] $H^{\ast}(G/T;\mathbb{R})=\mathbb{R}[\omega_{1},\cdots
,\omega_{n}]/\left\langle h_{i}^{(0)},g_{j}^{(0)};1\leq i\leq k,1\leq j\leq
m\right\rangle $,
\end{enumerate}

\noindent where $h_{i}^{(0)},g_{j}^{(0)}\in\mathbb{R}[\omega_{1},\cdots
,\omega_{n}]$ are the polynomials obtained from $h_{i}$, $g_{j}$ in Theorem
1.2 by substituting $-\frac{1}{p_{j}}\alpha_{j}$ for $y_{j}$, $1\leq j\leq m$.
For $G\neq E_{8}$ the result is shown by $n=k+m$ and (2.13). When $G=E_{8}$
the formula (2.13), together with the following relations by (1.5)

\begin{quote}
$g_{4}^{(0)}=\frac{6}{5}g_{6}^{(0)}$, $g_{7}^{(0)}=-\frac{3}{2}g_{6}^{(0)}$,
\end{quote}

\noindent implies that

\begin{enumerate}
\item[(2.14)] $H^{\ast}(E_{8}/T;\mathbb{R})=\mathbb{R}[\omega_{1}%
,\cdots,\omega_{8}]/\left\langle h_{i}^{(0)},g_{j}^{(0)};1\leq i\leq
k,j\in\{1,2,3,5,6\}\right\rangle $.
\end{enumerate}

\noindent The proof for $G=E_{8}$ is completed by $k=3$.

For the cases $\mathbb{F}=\mathbb{F}_{p}$ we set $G(p)=\{j\mid1\leq j\leq m$,
$p_{j}=p\}$ and let $\overline{G}(p)$ be the complement of $G(p)$ in
$\{1,\cdots,m\}$. After reduction $\operatorname{mod}p$ the relations of the
type $f_{j}$ in Theorem 1.2 turns to be

\begin{quote}
a) $\alpha_{j}\equiv0$ $\operatorname{mod}p$ if $j\in G(p)$;\quad

b) $y_{j}-q_{j}\alpha_{j}\equiv0$ $\operatorname{mod}p$ if $j\in\overline
{G}(p)$,
\end{quote}

\noindent where $q_{j}>0$ is the smallest integer satisfying $q_{j}p_{j}%
\equiv-1\operatorname{mod}p$. Formula a) implies that the relations $f_{j}$
with $j\in G(p)$ can be replaced by $\alpha_{j}\equiv0$. In view of b) we can
eliminate all the Schubert classes $y_{s}$ with $s\in\overline{G}(p)$ from the
set of generators to obtain the presentation

\begin{enumerate}
\item[(2.15)] $H^{\ast}(G/T;\mathbb{F}_{p})=\frac{\mathbb{F}_{p}[\omega
_{1},\cdots,\omega_{n},y_{j}]}{\left\langle h_{i}^{(p)},\alpha_{j}^{(p)}%
,g_{j}^{(p)},g_{s}^{(p)}\right\rangle _{1\leq i\leq k,j\in G(p),s\in
\overline{G}(p)}}$\textsl{ }with

i) $g_{j}^{(p)}=y_{j}^{k_{j}}+\beta_{j}^{(p)}$\textsl{ }with\textsl{ }%
$\beta_{j}^{(p)}\in\left\langle \omega_{1},\cdots,\omega_{n}\right\rangle
$,\textsl{ }$j\in G(p)$;\textsl{ }

ii) $\{h_{i}^{(p)},\alpha_{j}^{(p)},g_{s}^{(p)}\}\subset\left\langle
\omega_{1},\cdots,\omega_{n}\right\rangle $,
\end{enumerate}

\noindent where $h_{i}^{(p)},\alpha_{j}^{(p)},g_{s}^{(p)}$ and $\beta
_{j}^{(p)}$ are the polynomials obtained from $h_{i},\alpha_{j},g_{s}$ and
$\beta_{j}$, respectively, by substituting $q_{s}\alpha_{s}$ for $y_{s}$,
$s\in\overline{G}(p)$. For $G\neq E_{8}$ the lemma is verified by (2.15), as
well as the equality $n=k+\left\vert G(p)\right\vert +\left\vert \overline
{G}(p)\right\vert $.

For $G=E_{8}$ reduction $\operatorname{mod}p$ of the system (1.5) yields the relations

\begin{quote}
$g_{4}^{(p)}\equiv0$; $g_{6}^{(p)}\equiv y_{7}\alpha_{7}^{(p)}$
$\operatorname{mod}2$;

$g_{4}^{(p)}\equiv g_{7}^{(p)}\equiv-y_{6}^{2}\alpha_{6}^{(p)}$
$\operatorname{mod}3$;

$g_{6}^{(p)}\equiv g_{7}^{(p)}\equiv-y_{4}^{4}\alpha_{4}^{(p)}$
$\operatorname{mod}5$;

$g_{4}^{(p)}\equiv sg_{6}^{(p)}$; $g_{7}^{(p)}\equiv tg_{6}^{(p)}$
$\operatorname{mod}$ $p\neq2,3,5$ (for some $s,t\in\mathbb{F}_{p}$)
\end{quote}

\noindent among the system in ii). Therefore (2.15) implies that

\begin{enumerate}
\item[(2.16)] $H^{\ast}(E_{8}/T;\mathbb{F}_{p})=\frac{\mathbb{F}_{p}%
[\omega_{1},\cdots,\omega_{8},y_{j}]}{\left\langle e_{i}^{(p)},\alpha
_{j}^{(p)},g_{s}^{(p)},g_{j}^{(p)}\right\rangle _{1\leq i\leq3,j\in
E_{8}(p),s\in\widetilde{E}_{8}(p)}}$,
\end{enumerate}

\noindent where $\widetilde{E}_{8}(p)\subset\overline{E}_{8}(p)$ is the subset
defined by the table below

\begin{quote}%
\begin{tabular}
[c]{l|l|l|l|l}\hline\hline
$p$ & $2$ & $3$ & $5$ & $p\neq2,3,5$\\\hline
$\widetilde{E}_{8}(p)$ & $\overline{E}_{8}(2)\backslash\{4,6\}$ &
$\overline{E}_{8}(3)\backslash\{4,7\}$ & $\overline{E}_{8}(5)\backslash
\{6,7\}$ & $\overline{E}_{8}(p)\backslash\{4,7\}$\\\hline\hline
\end{tabular}
.
\end{quote}

\noindent The proof for the group $E_{8}$ is completed by $3+\left\vert
E_{8}(p)\right\vert +\left\vert \widetilde{E}_{8}(p)\right\vert =n$.$\square$

\bigskip

In term of the ring $H^{\ast}(G/T;\mathbb{F})$ presented through (2.13) to
(2.16) we call the set of the polynomials occurring in the relations that
belong also to the ideal $\left\langle \omega_{1},\cdots,\omega_{n}%
\right\rangle $ \textsl{the set of primary polynomials of }$G$\textsl{ over
}$\mathbb{F}$.

\bigskip

\noindent\textbf{Corollary 2.6.} \textsl{Let }$G$\textsl{ be a }%
$1$\textsl{--connected simple Lie group with rank }$n$\textsl{.}

\begin{quote}
\textsl{i) The set of primary polynomials of }$G$\textsl{ over }$\mathbb{R}%
$\textsl{ is}

$S_{0}(G)=\left\{
\begin{tabular}
[c]{l}%
$\{h_{i}^{(0)},g_{j}^{(0)}\mid1\leq i\leq k,1\leq j\leq m\}\text{ if }G\neq
E_{8}$;\\
$\{h_{i}^{(0)},g_{j}^{(0)}\mid1\leq i\leq3,j\in\{1,2,3,5,6\}\}\text{ if
}G=E_{8}\text{;}$%
\end{tabular}
\ \right.  $

\textsl{ii) The set of primary polynomials of }$G$\textsl{ over }%
$\mathbb{F}_{p}$\textsl{ is}

$S_{p}(G)=\left\{
\begin{tabular}
[c]{l}%
$\{h_{i}^{(p)},\alpha_{t}^{(p)},g_{s}^{(p)}\mid1\leq i\leq k,t\in
G(p),s\in\overline{G}(p)\}\text{ if }G\neq E_{8}$;\\
$\{h_{i}^{(p)},\alpha_{t}^{(p)},g_{s}^{(p)}\mid1\leq i\leq3,t\in E_{8}%
(p),s\in\widetilde{E_{8}}(p)\}\text{ if }G=E_{8}$.
\end{tabular}
\right.  $
\end{quote}

\noindent\textsl{In particular, in addition to }$\left\vert S_{0}%
(G)\right\vert =\left\vert S_{p}(G)\right\vert =n$\textsl{ one has}

\begin{enumerate}
\item[(2.17)] $\dim G=\sum\limits_{\alpha\in S_{0}(G)}(\left\vert
\alpha\right\vert -1)=\sum\limits_{\alpha\in S_{p}(G)}(\left\vert
\alpha\right\vert -1)+\sum\limits_{j\in G(p)}(k_{j}-1)\left\vert
y_{j}\right\vert .$
\end{enumerate}

\noindent\textbf{Proof.} The presentations of the sets $S_{0}(G)$ and
$S_{p}(G)$, together with the equalities $\left\vert S_{0}(G)\right\vert
=\left\vert S_{p}(G)\right\vert =n$, are transparent by the proof of Theorem
2.5. Concerning the equality (2.17) we remark that the dimension of all the
$1$--connected simple Lie groups are well known to be

\begin{quote}%
\begin{tabular}
[c]{l|llllllll}\hline\hline
$G$ & $SU(n)$ & $Spin(n)$ & $\qquad Sp(n)$ & $G_{2}$ & $F_{4}$ & $E_{6}$ &
$E_{7}$ & $E_{8}$\\\hline
$\dim G$ & $n^{2}-1$ & $\frac{n(n-1)}{2}$ & $\ n(2n+1)$ & $\ 14$ & $52$ & $78$
& $133$ & $248$\\\hline\hline
\end{tabular}
.
\end{quote}

\noindent On the other hand, one reads the degrees $\left\vert \alpha
\right\vert $ for $\alpha\in S_{0}(G)$, $S_{p}(G)$, as well as the data
$(k_{j},\left\vert y_{j}\right\vert )$ with $t\in G(p)$ of each simple Lie
group $G$ from the basic data given in Tables 5.1 and 5.2. When $\mathbb{F}%
=\mathbb{R}$ the equality (2.17) is obvious. For $\mathbb{F}=\mathbb{F}_{p}$
it can be verified in accordance to $p\in$ $\left\{  2,3,5\right\}  $ and
$p\notin\left\{  2,3,5\right\}  $.$\square$

\subsection{The Bockstein cohomology of\textsl{ }a space $X$}

\noindent Recall that any finitely generated abelian group $A$ admits a decomposition

\begin{enumerate}
\item[(2.18)] $A=\mathcal{F}(A)\underset{p}{\oplus}\sigma_{p}(A)$ with
$\sigma_{p}(A):=\{x\in A\mid p^{r}x=0,$ $r\geq1\}$,
\end{enumerate}

\noindent where the sum is over all primes $p>1$, and where the components
$\sigma_{p}(A)$ and $\mathcal{F}(A)$ are called, respectively, the
$p$--\textsl{primary component} and a \textsl{free part }of the group $A$. If
$A$ is the integral cohomology $H^{\ast}(X)$ of a space $X$ we use
$\mathcal{F}(X)$ in place of $\mathcal{F}(H^{\ast}(X))$, and write $\sigma
_{p}(X)$ instead of $\sigma_{p}(H^{\ast}(X))$.

Given a topological space $X$ and a prime $p$ the \textsl{Bockstein operator}
$\delta_{p}$ on the cohomology $H^{\ast}(X;\mathbb{F}_{p})$ is the
differential with degree $1$

\begin{quote}
$\delta_{p}:=r_{p}\circ\beta_{p}:H^{r}(X;\mathbb{F}_{p})\overset{\beta_{p}%
}{\rightarrow}H^{r+1}(X)\overset{r_{p}}{\rightarrow}H^{r+1}(X;\mathbb{F}_{p})$,
\end{quote}

\noindent where $r_{p}$ is the module $p$ reduction and $\beta_{p}$ is the
Bockstein homomorphism. The pair $\{H^{\ast}(X;\mathbb{F}_{p});\delta_{p}\}$
is a cochain complex whose cohomology, denoted by $\overline{H}^{\ast
}(X;\mathbb{F}_{p})$, is called the $\operatorname{mod}p$ \textsl{Bockstein
cohomology }of the space\textsl{ }$X$.

\bigskip

\noindent\textbf{Theorem 2.7. }\textsl{For a topological space }%
$X$\textsl{\ and a prime }$p$\textsl{\ one has}

\begin{enumerate}
\item[(2.19)] $\dim_{\mathbb{R}}\mathcal{F}(X)\otimes\mathbb{R}\leq
\dim_{\mathbb{F}_{p}}\overline{H}^{\ast}(X;\mathbb{F}_{p})$,
\end{enumerate}

\noindent\textsl{while the equality implies that}

\begin{quote}
\textsl{i)} $\sigma_{p}(X)=\operatorname{Im}\beta_{p}$;\quad\textsl{ii)
}$\sigma_{p}(X)\overset{\cong}{\rightarrow}\operatorname{Im}\delta_{p}$
\textsl{via} $r_{p}$\textsl{.}
\end{quote}

\noindent\textbf{Proof. }Let $r_{p}^{1}$ be the restriction of $r_{p}$ on
$\operatorname{Im}\beta_{p}$, and consider the short exact ladder induced by
the decomposition $\delta_{p}=r_{p}\circ\beta_{p}$

\begin{quote}%
\begin{tabular}
[c]{lllllll}%
$0\rightarrow$ & $\operatorname{Im}r_{p}$ & $\rightarrow$ & $H^{\ast
}(X;\mathbb{F}_{p})$ & $\overset{\beta_{p}}{\rightarrow}$ & $\operatorname{Im}%
\beta_{p}$ & $\rightarrow0$\\
& $\quad i\downarrow$ &  & $\quad\parallel$ &  & $r_{p}^{1}\downarrow$ & \\
$0\rightarrow$ & $\ker\delta_{p}$ & $\rightarrow$ & $H^{\ast}(X;\mathbb{F}%
_{p})$ & $\overset{\delta_{p}}{\rightarrow}$ & $\operatorname{Im}\delta_{p}$ &
$\rightarrow0$%
\end{tabular}
.
\end{quote}

\noindent Since the vertical map $i$ on the left is injective while the
vertical map $r_{p}^{1}$ on the right is surjective, one has by \cite[Exercise
5, p.726]{Wh} that

\begin{quote}
$\ker r_{p}^{1}\cong co\ker i$.
\end{quote}

\noindent In addition, with respect to the obvious decompositions

\begin{quote}
$\operatorname{Im}r_{p}\cong\mathcal{F}(X)/p\cdot\mathcal{F}(X)\oplus
\sigma_{p}(X)/p\cdot\sigma_{p}(X)$,

$\ker\delta_{p}=\overline{H}^{\ast}(X;\mathbb{F}_{p})\oplus\operatorname{Im}%
\delta_{p}$
\end{quote}

\noindent the monomorphism $i$ satisfies the the following relations

\begin{quote}
a) $i(\mathcal{F}(X)/p\cdot\mathcal{F}(X))\subseteq\overline{H}^{\ast
}(X;\mathbb{F}_{p})$;

b) $\operatorname{Im}\delta_{p}\subseteq i(\sigma_{p}(X)/p\cdot\sigma_{p}(X))$
(since $\operatorname{Im}\beta_{p}\subseteq\sigma_{p}(X)$)
\end{quote}

\noindent With $\dim_{\mathbb{F}_{p}}\frac{\mathcal{F}(X)}{p\cdot
\mathcal{F}(X)}=\dim_{\mathbb{R}}\mathcal{F}(X)\otimes\mathbb{R}$ formula
(2.19) is shown by a).

If the equality in (2.19) (i.e. in a)) holds then $i$ is surjective by a) and
b). Hence $r_{p}^{1}$ is an isomorphism. We get properties i) and ii) from
$\operatorname{Im}\beta_{p}\subseteq\sigma_{p}(X)$, as well as the fact that
$\ker r_{p}^{1}=0$ implies that $\operatorname{Im}\beta_{p}=\sigma_{p}%
(X)$.$\square$

\subsection{Canonical maps $E_{3}^{\ast,0}(G;\mathcal{R})$, $E_{3}^{\ast
,1}(G;\mathcal{R})\rightarrow H^{\ast}(G;\mathcal{R})$}

The spectral sequence $\{E_{r}^{\ast,\ast}(G;\mathcal{R}),d_{r}\}$ has certain
crucial properties that allow us to solve the extension problem from
$E_{3}^{\ast,\ast}(G;\mathcal{R})$ to $H^{\ast}(G;\mathcal{R})$. To be precise
let $F^{p}$ be the filtration on $H^{\ast}(G;\mathcal{R})$ defined by (1.1).
That is

\begin{center}
$0=F^{r+1}(H^{r}(G;\mathcal{R}))\subseteq F^{r}(H^{r}(G;\mathcal{R}%
))\subseteq\cdots\subseteq F^{0}(H^{r}(G;\mathcal{R}))=H^{r}(G;\mathcal{R})$
\end{center}

\noindent with $E_{\infty}^{p,q}(G;\mathcal{R})=F^{p}(H^{p+q}(G;\mathcal{R}%
))/F^{p+1}(H^{p+q}(G;\mathcal{R}))$.

The routine relation $d_{r}(E_{r}^{\ast,0}(G;\mathcal{R}))=0$ for $r\geq2$
yields the sequence of quotient maps

\begin{center}
$H^{r}(G/T;\mathcal{R})=E_{2}^{r,0}\rightarrow E_{3}^{r,0}\rightarrow
\cdots\rightarrow E_{\infty}^{r,0}=F^{r}(H^{r}(G;\mathcal{R}))\subset
H^{r}(G;\mathcal{R})$
\end{center}

\noindent whose composition agrees with the induces map $\pi^{\ast}:H^{\ast
}(G/T;\mathcal{R})\rightarrow H^{\ast}(G;\mathcal{R})$ \cite[P.147]{Mc}. For
this reason we can reserve $\pi^{\ast}$ for the composition

\begin{enumerate}
\item[(2.20)] $\pi^{\ast}:E_{3}^{\ast,0}(G;\mathcal{R})\rightarrow
\cdots\rightarrow E_{\infty}^{\ast,0}(G;\mathcal{R})=F^{r}(H^{r}%
(G;\mathcal{R}))\subset H^{\ast}(G;\mathcal{R})$.
\end{enumerate}

The property $H^{odd}(G/T;\mathcal{R})=0$ by Theorem 1.2 indicates the relations

\begin{quote}
a) $E_{r}^{odd,q}=0$ for $r,q\geq0$;\quad

b) $E_{3}^{4k,2}=E_{4}^{4k,2}=\cdots=E_{\infty}^{4k,2}$.
\end{quote}

\noindent From $F^{2k+1}(H^{2k+1}(G;\mathcal{R}))=F^{2k+2}(H^{2k+1}%
(G;\mathcal{R}))=0$ by a) one finds that

\begin{quote}
$E_{\infty}^{2k,1}(G;\mathcal{R})=F^{2k}(H^{2k+1}(G;\mathcal{R}))\subset
H^{2k+1}(G;\mathcal{R})$.
\end{quote}

\noindent Combining this with $d_{r}(E_{r}^{\ast,1})=0$ for $r\geq3$ yields
the composition

\begin{enumerate}
\item[(2.21)] $\kappa:E_{3}^{\ast,1}(G;\mathcal{R})\rightarrow\cdots
\rightarrow E_{\infty}^{\ast,1}(G;\mathcal{R})\subset H^{2t+1}(G;\mathcal{R})$
\end{enumerate}

\noindent that interprets elements of $E_{3}^{\ast,1}$ as cohomology classes
of the group $G$.

Assume that $\dim G/T=g$ and that the rank of $G$ is $n$. Since $E_{2}%
^{s,t}(G;\mathcal{R})=0$ for $s>g$ or $t>n$ by (1.2), any differential $d_{r}$
that acts or lands on the group $E_{r}^{g,n}(G;\mathcal{R})$ must be trivial.
As a result one gets

\begin{enumerate}
\item[(2.22)] $E_{r}^{g,n}(G;\mathcal{R})=H^{\dim G}(G;\mathcal{R}%
)=\mathcal{R}$ for all $r\geq2$.
\end{enumerate}

Finally, with the multiplicative structure inherited from that on $E_{2}%
^{\ast,\ast}(G;\mathcal{R})$ the third page $E_{3}^{\ast,\ast}(G;\mathcal{R})$
is a bi-graded ring \cite[P.668]{Wh}. Granted with the maps $\pi^{\ast}$ and
$\kappa$ in (2.20) and (2.21) we show that

\bigskip

\noindent\textbf{Lemma 2.8. }\textsl{For any }$\xi\in\operatorname{Im}\kappa
$\textsl{\ one has }$\xi^{2}\in\operatorname{Im}\pi^{\ast}$\textsl{.}

\bigskip

\noindent\textbf{Proof.} For an element $x\in E_{3}^{2k,1}$ the obvious
relation $x^{2}=0$ in

\begin{quote}
$E_{3}^{4k,2}=E_{\infty}^{4k,2}=\mathcal{F}^{4k}H^{4k+2}/\mathcal{F}%
^{4k+1}H^{4k+2}$ (by b)).
\end{quote}

\noindent implies that $\kappa(x)^{2}\in\mathcal{F}^{4k+1}H^{4k+2}$. From

\begin{quote}
$\mathcal{F}^{4k+1}H^{4k+2}/\mathcal{F}^{4k+2}H^{4k+2}=E_{\infty}^{4k+1,1}=0 $
(by a))
\end{quote}

\noindent one gets further that $\kappa(x)^{2}\in\mathcal{F}^{4k+2}H^{4k+2}$.
The proof is completed by

\begin{quote}
$\mathcal{F}^{4k+2}H^{4k+2}=E_{\infty}^{4k+2,0}=\operatorname{Im}\pi^{\ast}%
$.$\square$
\end{quote}

\bigskip

We conclude this section with a standard result of algebra. Let $A=\oplus
_{i\geq0}A^{i}$ be a graded ring over $\mathcal{R}$, and let $u=t_{1}^{a_{1}%
}\cdots t_{k}^{a_{k}}\in A$ be a decomposed element with $a_{i}\geq1$. The
element $u$ is called \textsl{monotonous} in $A$ if the following two
conditions are met

\begin{quote}
i) the group $A^{r}=\mathcal{R}$ is spanned by $u$, where $r=\left\vert
u\right\vert $;

ii) $t_{1}^{b_{1}}\cdots t_{k}^{b_{k}}=0$ if $\left\vert t_{1}^{b_{1}}\cdots
t_{k}^{b_{k}}\right\vert =r$ but $(b_{1},\cdots,b_{k})\neq(a_{1},\cdots
,a_{k})$.
\end{quote}

\noindent\textbf{Lemma 2.9.} \textsl{If }$u=t_{1}^{a_{1}}\cdots t_{k}^{a_{k}%
}\in A$ \textsl{is a monotonous element, then the set }$\{t_{1}^{r_{1}}\cdots
t_{k}^{r_{k}}\mid0\leq r_{i}\leq a_{i}\}$\textsl{ of monomials is linearly
independent in }$A$\textsl{, and spans a direct summand of the ring }%
$A$\textsl{.}$\square$

\section{The additive structure on $H^{\ast}(G;\mathcal{R})$}

Assume that $G$ is a $1$--connected simple Lie group with rank $n$. As in
(2.21) we set $g:=\dim G/T$. In Theorems 3.4 and 3.6 we formulate the additive
cohomology $H^{\ast}(G;\mathcal{R})$ in term of the presentation (1.4) of the
ring $H^{\ast}(G/T)$.

\subsection{The ring $E_{3}^{\ast,0}(G;\mathcal{R})$}

By (1.2) and (1.3) the group $E_{3}^{\ast,0}(G;\mathcal{R})$ is the cokernel
of the differential

\begin{quote}
$d_{2}:H^{\ast}(G/T;\mathcal{R})\otimes\Lambda^{1}(t_{1},\cdots,t_{n}%
)\rightarrow H^{\ast}(G/T;\mathcal{R})$.
\end{quote}

\noindent With $\operatorname{Im}d_{2}=$ $\left\langle \omega_{1}%
,\cdots,\omega_{n}\right\rangle $ we get from Theorem 1.2, together with the
presentations of the rings $H^{\ast}(G/T;\mathbb{F})$ in (2.13) and (2.15), that

\bigskip

\noindent\textbf{Lemma 3.1.} \textsl{Let }$\left\{  y_{1},\cdots
,y_{m}\right\}  $\textsl{\ be the special Schubert classes on }$G/T$
\textsl{given by Theorem 1.2. Then}

\begin{enumerate}
\item[(3.1)] $E_{3}^{\ast,0}(G;\mathcal{R})=\left\{
\begin{tabular}
[c]{l}%
$\mathbb{Z}[y_{1},\cdots,y_{m}]/\left\langle p_{j}y_{j},y_{j}^{k_{j}%
}\right\rangle _{1\leq j\leq m}$ \textsl{if} $\mathcal{R}=\mathbb{Z}%
$\textsl{;}\\
$\mathbb{F}_{p}[y_{j}]_{j\in G(p)}/\left\langle y_{j}^{k_{j}}\right\rangle $
\textsl{if} $\mathcal{R}=\mathbb{F}_{p}$\textsl{;}\\
$\mathbb{R}$\textsl{ if }$\mathcal{R}=\mathbb{R}$\textsl{.}$\square$%
\end{tabular}
\right.  \ \ $
\end{enumerate}

\noindent\textbf{Example 3.2. }Let $G$ be an exceptional Lie group. Inputting
the values of the integers $p_{i}$, $k_{i}$ given by Table 5.2 into (3.1)
yields the following formulae of the ring $E_{3}^{\ast,0}(G)$, where to
emphasize the degrees of the generators the notion $x_{\left\vert
y_{i}\right\vert }$ is used in place of $y_{i}$.

\begin{quote}
$E_{3}^{\ast,0}(G_{2})=\mathbb{Z}[x_{6}]/\left\langle 2x_{6},x_{6}%
^{2}\right\rangle $;

$E_{3}^{\ast,0}(F_{4})=\mathbb{Z}[x_{6},x_{8}]/\left\langle 2x_{6},x_{6}%
^{2},3x_{8},x_{8}^{3}\right\rangle $;

$E_{3}^{\ast,0}(E_{6})=\mathbb{Z}[x_{6},x_{8}]/\left\langle 2x_{6},x_{6}%
^{2},3x_{8},x_{8}^{3}\right\rangle $;

$E_{3}^{\ast,0}(E_{7})=\mathbb{Z}[x_{6},x_{8},x_{10},x_{18}]/\left\langle
2x_{6},3x_{8},2x_{10},2x_{18},x_{6}^{2},x_{8}^{3},x_{10}^{2},x_{18}%
^{2}\right\rangle $;

$E_{3}^{\ast,0}(E_{8})=\frac{\mathbb{Z}[x_{6},x_{8},x_{10},x_{12}%
,x_{18},x_{20},x_{30}]}{\quad\left\langle 2x_{6},3x_{8},2x_{10},5x_{12}%
,2x_{18},3x_{20},2x_{30},x_{6}^{8},x_{8}^{3},x_{10}^{4},x_{12}^{5},x_{18}%
^{2},x_{20}^{3},x_{30}^{2}\right\rangle }$.$\square$
\end{quote}

\subsection{The additive cohomology $H^{\ast}(G;\mathbb{F})$ with
$\mathbb{F}=\mathbb{R}$ or $\mathbb{F}_{p}$}

Let $S_{p}(G)$ ($S_{0}(G)$) be the set of primary polynomials of $G$ over
$\mathbb{F}_{p}$ (over $\mathbb{R}$), see Corollary 2.6. In term of the map
$\varphi$ in Lemma 2.2 we set

\begin{enumerate}
\item[(3.1)] $\xi_{\left\vert \alpha\right\vert -1}:=\varphi(\alpha)\in
E_{3}^{\left\vert \alpha\right\vert -2,1}(G;\mathbb{F})$, $\alpha\in S_{p}(G)$
(resp. $\alpha\in S_{0}(G)$).
\end{enumerate}

\noindent This notion is justified by the fact that $\left\vert \alpha
\right\vert \neq\left\vert \beta\right\vert $ for any pair $\alpha,\beta\in
S_{p}(G)$ with $\alpha\neq\beta$, see in Tables 5.1 and 5.2. By Lemma 2.5 and
Theorem 2.4 the inclusions $\xi_{\left\vert \alpha\right\vert -1}\in
E_{3}^{\ast,\ast}(G;\mathbb{F})$, $y_{t}\in E_{3}^{\ast,\ast}(G;\mathbb{F}%
_{p})$ with $t\in G(p)$, extend to a surjective ring map

\begin{enumerate}
\item[(3.2)] $f_{p}:\mathbb{F}_{p}[y_{t}]_{t\in G(p)}/\left\langle
y_{t}^{k_{t}}\right\rangle \otimes\Lambda_{\mathbb{F}_{p}}(\xi_{\left\vert
\alpha\right\vert -1})_{\alpha\in S_{p}(G)}\rightarrow E_{3}^{\ast,\ast
}(G;\mathbb{F}_{p})$

(resp. $f_{0}:\Lambda_{\mathbb{R}}(\xi_{\left\vert \alpha\right\vert
-1})_{\alpha\in S_{0}(G)}\rightarrow E_{3}^{\ast,\ast}(G;\mathbb{R})$).
\end{enumerate}

\noindent\textbf{Lemma 3.3. }\textsl{The map }$f_{p}$\textsl{ (resp. }$f_{0}%
$\textsl{) in (3.2) is a ring isomorphism.}

\bigskip

\noindent\textbf{Proof.} It suffices to show that the map $f_{p}$ ($f_{0}$) is
also injective. Consider firstly the map $f_{0}$. In the top degree the ring
$\Lambda_{\mathbb{R}}(\xi_{\left\vert \alpha\right\vert -1})_{\alpha\in
S_{0}(G)}\mathbb{\ }$is generated by the single element $u=\prod_{\alpha\in
S_{0}(G)}\xi_{\left\vert \alpha\right\vert -1}$. Since $\left\vert
u\right\vert =\dim G(=g+n)$ by (2.17) the class $f_{0}(u)\in E_{3}%
^{g,n}(G;\mathbb{R})=\mathbb{R}$ must be a basis element by the surjectivity
of $f_{0}$. With $\left\vert S_{0}(G)\right\vert =n$ by Corollary 2.6 the
result is shown by

\begin{quote}
$2^{n}=\dim\Lambda_{\mathbb{R}}(\xi_{\left\vert \alpha\right\vert -1}%
)_{\alpha\in S_{0}(G)}\geq\dim E_{3}^{\ast,\ast}(G;\mathbb{R})$ $\geq2^{n}$,
\end{quote}

\noindent where the first inequality comes from the surjectivity of $f_{0}$,
the second is obtained by Lemma 2.9, since the class $f_{0}(u)=$
$\prod_{\alpha\in S_{0}(G)}f_{0}(\xi_{\left\vert \alpha\right\vert -1})\in$
$E_{3}^{\ast,\ast}(G;\mathbb{R})$ is monotonous\textsl{ }by (2.22) and Lemma 2.8.

The same argument applies equally well to the map $f_{p}$. In the top degree
the ring $\mathbb{F}_{p}[y_{t}]_{t\in G(p)}/\left\langle y_{t}^{k_{t}%
}\right\rangle \otimes\Lambda_{\mathbb{F}_{p}}(\xi_{\left\vert \alpha
\right\vert -1})_{\alpha\in S_{p}(G)}$ has the basis element

\begin{enumerate}
\item[(3.3)] $u_{p}=\prod\limits_{t\in G(p)}y_{t}^{k_{t}-1}\prod
\limits_{\alpha\in S_{p}(G)}\xi_{\left\vert \alpha\right\vert -1}$.
\end{enumerate}

\noindent Since $\deg u_{p}=\dim G(=g+n)$ by (2.17) the class $f_{p}(u_{p})\in
E_{3}^{g,n}(G;\mathbb{F}_{p})$ $=\mathbb{F}_{p}$ must be a basis element. The
proof is done by the relations

\begin{quote}
$\dim\mathbb{F}_{p}[y_{t}]_{t\in G(p)}/\left\langle y_{t}^{k_{t}}\right\rangle
\otimes\Lambda_{\mathbb{F}_{p}}(\xi_{\left\vert \alpha\right\vert -1}%
)_{\alpha\in S_{p}(G)}=2^{n}\prod\limits_{t\in G(p)}k_{t}$

$\quad\geq\dim E_{3}^{\ast,\ast}(G;\mathbb{F}_{p})\geq2^{n}\prod\limits_{t\in
G(p)}k_{t}$,
\end{quote}

\noindent where the first inequality comes from the surjectivity of $f_{p}$,
while the second is obtained by Lemma 2.9, since the element $f_{p}(u_{p})\in
E_{3}^{\ast,\ast}(G;\mathbb{F}_{p})$ is monotonous\textsl{ }by (2.22) and
Lemma 2.8.$\square$

\bigskip

Since the third page $E_{3}^{\ast,\ast}(G;\mathbb{F})$ is generated
multiplicatively by $E_{3}^{\ast,0}(G;\mathbb{F})$ and $E_{3}^{\ast
,1}(G;\mathbb{F})$ by Lemma 3.3 we have $E_{\infty}^{\ast,\ast}(G;\mathbb{F}%
)=E_{3}^{\ast,\ast}(G;\mathbb{F})$. In particular, the maps $\pi^{\ast}$ and
$\kappa$ in (2.20) and (2.21) are monomorphisms of $\mathbb{F}$-spaces.
Combining this with the isomorphism $E_{\infty}^{\ast,\ast}(G;\mathbb{F}%
)=H^{\ast}(G;\mathbb{F})$ of $\mathbb{F}$--spaces we get from Lemma 3.3 the
following additive presentation of the cohomology $H^{\ast}(G;\mathbb{F})$,
where for simplicity the notion $\xi_{\left\vert \alpha\right\vert -1}$ and
$x_{\left\vert y_{t}\right\vert }$ are reserved for the cohomology classes
$\kappa(\xi_{\left\vert \alpha\right\vert -1})\in H^{\ast}(G;\mathbb{F}_{p})$
and $\pi^{\ast}(y_{t})\in H^{\ast}(G;\mathbb{F}_{p})$ (see (2.20) and (2.21)),
respectively, and where we need to replace exterior ring $\Lambda
_{\mathbb{F}_{p}}(\xi_{\left\vert \alpha\right\vert -1})$ by the
$\mathbb{F}_{p}$--module $\Delta_{\mathbb{F}_{p}}(\xi_{\left\vert
\alpha\right\vert -1})$ since the property $\xi_{\left\vert \alpha\right\vert
-1}^{2}=0$ on $E_{\infty}^{\ast,\ast}(G;\mathbb{F}_{p})$ may not survive to
$H^{\ast}(G;\mathbb{F}_{p})$, see Lemma 2.8.

\bigskip

\noindent\textbf{Theorem 3.4.} \textsl{Let }$G$\textsl{ be an }$1$%
\textsl{--connected simple Lie group with} \textsl{the set }$S_{p}(G)$
\textsl{(resp. }$S_{0}(G)$\textsl{)} \textsl{of primary polynomials. Then}

\begin{enumerate}
\item[(3.4)] $H^{\ast}(G;\mathbb{F})=\left\{
\begin{tabular}
[c]{l}%
$\mathbb{F}_{p}[x_{\left\vert y_{t}\right\vert }]_{t\in G(p)}/\left\langle
x_{\left\vert y_{t}\right\vert }^{k_{t}}\right\rangle \otimes\Delta
_{\mathbb{F}_{p}}(\xi_{\left\vert \alpha\right\vert -1})_{\alpha\in S_{p}(G)}$
\textsl{if} $\mathbb{F}=\mathbb{F}_{p}$\textsl{;}\\
$\Delta_{\mathbb{R}}(\xi_{\left\vert \alpha\right\vert -1})_{\alpha\in
S_{0}(G)}$ \textsl{if} $\mathbb{F}=\mathbb{R}$.$\square$%
\end{tabular}
\right.  \ \ \ \ \ \ $
\end{enumerate}

\subsection{The torsion ideal $\sigma_{p}(G)\subset H^{\ast}(G)$}

To analyze the Bockstein operator on $H^{\ast}(G;\mathbb{F}_{p})$

\begin{quote}
$\delta_{p}:=r_{p}\circ\beta_{p}:H^{r}(G;\mathbb{F}_{p})\overset{\beta_{p}%
}{\rightarrow}H^{r+1}(G)\overset{r_{p}}{\rightarrow}H^{r+1}(G;\mathbb{F}_{p})$.
\end{quote}

\noindent with respect to the presentation (3.4) introduce the partition on
$S_{p}(G)$ by

\begin{quote}
$S_{p}(G)=S_{p}^{1}(G)\amalg S_{p}^{2}(G)$ with $S_{p}^{1}(G)=\{\alpha
_{t}^{(p)}\mid t\in G(p)\}$,
\end{quote}

\noindent see Corollary 2.6. For each subset $I\subset S_{p}^{1}(G)$ define
the following elements of $H^{\ast}(G;\mathbb{F}_{p})$

\begin{enumerate}
\item[(3.5)] $\xi_{I}=\prod\limits_{\alpha\in I}\xi_{\left\vert \alpha
\right\vert -1}$; $c_{I}=\delta_{p}\xi_{I}$, $D_{I}=\sum\limits_{t\in
I}-x_{\left\vert y_{t}\right\vert }c_{I_{t}}$, $R_{I}=(\prod\limits_{t\in
I}^{k_{t}-1}x_{\left\vert y_{t}\right\vert }^{k_{t}-1})c_{I}$.
\end{enumerate}

\noindent\textbf{Theorem 3.5.} \textsl{Let }$\sigma_{p}(G)$ \textsl{be the
}$p$\textsl{--primary component of the group} $H^{\ast}(G)$\textsl{.}

\begin{quote}
\textsl{i) The map }$r_{p}$ \textsl{restricts to an isomorphism }$\sigma
_{p}(G)\cong\operatorname{Im}\delta_{p}$;

\textsl{ii) The subgroup} $\operatorname{Im}\delta_{p}\subset H^{\ast
}(G;\mathbb{F}_{p})$ \textsl{has the presentation}

$\qquad\operatorname{Im}\delta_{p}=\frac{\mathbb{F}_{p}[x_{\left\vert
y_{t}\right\vert }]_{t\in G(p)}/\left\langle x_{\left\vert y_{t}\right\vert
}^{k_{t}}\right\rangle \{1,c_{I}\}^{+}}{\left\langle D_{J},R_{K}\right\rangle
}\otimes\Delta(\xi_{\left\vert \alpha\right\vert -1})_{\alpha\in S_{p}^{2}%
(G)}$\textsl{, }
\end{quote}

\noindent\textsl{where} $I,J,K\subseteq S_{p}^{1}(G)$\textsl{\ with
}$\left\vert I\right\vert ,\left\vert J\right\vert ,\left\vert K\right\vert
\geq2$\textsl{.}

\bigskip

\noindent\textbf{Proof. }As the ring $H^{\ast}(G/T)$ is torsion free one has
the short exact sequence of Koszul complexes

\begin{quote}
$0\rightarrow$ $E_{2}^{\ast,\ast}(G)\overset{\cdot p}{\rightarrow}E_{2}%
^{\ast,\ast}(G)\overset{r_{p}}{\rightarrow}E_{2}^{\ast,\ast}(G;\mathbb{F}%
_{p})\rightarrow0$.
\end{quote}

\noindent The corresponding connecting homomorphism $\widehat{\beta}_{p}%
:E_{3}^{\ast,1}(G;\mathbb{F}_{p})\rightarrow E_{3}^{\ast,0}(G)$ can be shown
to satisfy the following formulae

\begin{quote}
a) $\widehat{\beta}_{p}(\xi_{\left\vert \alpha\right\vert -1})=-x_{\left\vert
\alpha\right\vert }$ if $\alpha\in S_{p}^{1}(G)$; $0$ if $\alpha\in S_{p}%
^{2}(G)$.
\end{quote}

\noindent Indeed, for each $t\in G(p)$ the diagram chasing by the relation
$f_{t}$ on $H^{\ast}(G/T)$

\begin{quote}
$%
\begin{array}
[c]{ccccc}
&  & \widetilde{\alpha}_{t} & \rightarrow & \widetilde{\alpha}_{t}^{(p)}\\
&  & d_{2}\downarrow\quad &  & d_{2}\downarrow\quad\\
-y_{t} & \overset{p}{\longrightarrow} & \alpha_{t} &  & \alpha_{t}^{(p)}=0
\end{array}
$
\end{quote}

\noindent in above exact sequence shows formula a) for $\alpha\in S_{p}%
^{1}(G)$, while the formula for $\alpha\in S_{p}^{2}(G)$ comes from the relations

\begin{quote}
$r_{p}(h_{i})=h_{i}^{(p)}$, $1\leq i\leq k$; $r_{p}(g_{s})=g_{s}^{(p)}$,
$s\in\overline{G}(p)$
\end{quote}

\noindent by the proof of Theorem 2.5. Moreover, with respect to the
monomorphisms $\kappa$ and $\pi^{\ast}$ the map $\widehat{\beta}_{p}$ clearly
fits the commutative diagram

\begin{quote}%
\begin{tabular}
[c]{lll}%
$E_{3}^{\ast,1}(G;\mathbb{F}_{p})$ & $\overset{\widehat{\beta}_{p}%
}{\rightarrow}$ & $E_{3}^{\ast,0}(G)$\\
$\kappa\downarrow$ &  & $\pi^{\ast}\downarrow$\\
$H^{\ast}(G;\mathbb{F}_{p})$ & $\overset{\beta_{p}}{\rightarrow}$ & $H^{\ast
}(G)$%
\end{tabular}

\end{quote}

\noindent that allows us to translate a) into the formulae for $\beta_{p}$

\begin{quote}
b) $\beta_{p}(\xi_{\left\vert \alpha\right\vert -1})=-x_{\left\vert
\alpha\right\vert }$ if $\alpha\in S_{p}^{1}(G)$; $0$ if $\alpha\in S_{p}%
^{2}(G)$.
\end{quote}

\noindent It implies that if we set $C^{\ast}=\frac{\mathbb{F}_{p}%
[x_{\left\vert y_{t}\right\vert }]_{t\in G(p)}}{\left\langle x_{\left\vert
y_{t}\right\vert }^{k_{t}}\right\rangle }\otimes\Delta(\xi_{\left\vert
\alpha\right\vert -1})_{\alpha\in S_{p}^{1}(G)}$ and write

\begin{quote}
$H^{\ast}(G;\mathbb{F}_{p})=C^{\ast}\otimes\Delta(\xi_{\left\vert
\alpha\right\vert -1})_{\alpha\in S_{p}^{2}(G)}$ (by Theorem 3.4)
\end{quote}

\noindent then both factors $C^{\ast}$ and $\Delta(\xi_{\left\vert
\alpha\right\vert -1})_{\alpha\in S_{p}^{2}(G)}$ are closed subspaces of
$\delta_{p}$. In addition, formula b) indicates that

\begin{quote}
c) $\{C^{\ast},\delta_{p}\}$ is the Koszul complex concerned by Theorem 2.1;

d) $\delta_{p}$ acts trivially on the second factor $\Delta(\xi_{\left\vert
\alpha\right\vert -1})_{\alpha\in S_{p}^{2}(G)}$.
\end{quote}

\noindent Theorem 2.1, together with K\"{u}nneth formula, gives rise to the
formula in ii), as well as the decomposition on $\overline{H}^{\ast
}(G;\mathbb{F}_{p})$

\begin{quote}
$\overline{H}^{\ast}(G;\mathbb{F}_{p})=\Delta(x_{\left\vert \alpha\right\vert
}^{k_{t}-1}\xi_{\left\vert \alpha\right\vert -1})_{\alpha\in S_{p}^{1}%
(G)}\otimes\Delta(\xi_{\left\vert \alpha\right\vert -1})_{\alpha\in S_{p}%
^{2}(G)}$.
\end{quote}

\noindent With $\dim_{\mathbb{F}_{p}}\overline{H}^{\ast}(G;\mathbb{F}%
_{p})=2^{n}$($=\dim_{\mathbb{R}}H^{\ast}(G;\mathbb{R})$) by $\left\vert
S_{p}(G)\right\vert =n$ the map $r_{p}$ is an isomorphism by Theorem 2.7. This
completes the proof.$\square$

\subsection{The additive cohomology $H^{\ast}(G)$}

By (2.18) the integral cohomology $H^{\ast}(G)$ admits the decomposition

\begin{quote}
$H^{\ast}(G)=\mathcal{F}(G)\underset{p}{\oplus}\sigma_{p}(G)$.
\end{quote}

\noindent in which the primary components $\sigma_{p}(G)$ have been determined
by Theorem 3.5. To clarify the free part $\mathcal{F}(G)$ we introduce in term
of Theorem 1.2 \textsl{the set }$S(G)$\textsl{ of primary polynomials of }%
$G$\textsl{ over the ring} $\mathbb{Z}$ by

\begin{enumerate}
\item[(3.6)] $S(G):=\left\{  h_{i},p_{j}g_{j}-y_{j}^{k_{j}-1}f_{j}\mid1\leq
i\leq k\right\}  $,
\end{enumerate}

\noindent where $1\leq j\leq m$ with $j\neq4,7$ if $G=E_{8}$. Since

\begin{quote}
$S(G)\subset\left\langle \omega_{1},\cdots,\omega_{n}\right\rangle
\cap\left\langle h_{i},f_{j},g_{j}\right\rangle _{1\leq i\leq k;1\leq j\leq
m}$\textsl{ }with\textsl{ }$\left\vert S(G)\right\vert =n$
\end{quote}

\noindent the map $\varphi$ in Lemma 2.2 yields the set of $n$ elements in
$E_{3}^{\ast,1}(G)$

\begin{quote}
$\left\{  \rho_{\left\vert \alpha\right\vert -1}:=\varphi(\alpha)\in
E_{3}^{\ast,1}(G)\mid\alpha\in S(G)\right\}  $.
\end{quote}

\noindent Again, the notion $\rho_{\left\vert \alpha\right\vert -1}$ is valid
as the polynomials of the set $S(G)$ have distinct degrees. For simplicity we
reserve the notion $\rho_{\left\vert \alpha\right\vert -1}$ also for the
cohomology class $\kappa(\rho_{\left\vert \alpha\right\vert -1})\in H^{\ast
}(G)$. The following result assures us that the square free monomials in
$\rho_{\left\vert \alpha\right\vert -1}$ form a basis of $\mathcal{F}(G)$.

\bigskip

\noindent\textbf{Theorem 3.6.} \textsl{If }$G$\textsl{ is a }$1$%
\textsl{--connected simple Lie group, then}

\begin{enumerate}
\item[(3.7)] $H^{\ast}(G)=\Delta_{\mathbb{Z}}(\rho_{\left\vert \alpha
\right\vert -1})_{\alpha\in S(G)}\underset{p\in\{2,3,5\}}{\oplus}\sigma
_{p}(G)$\textsl{,}
\end{enumerate}

\noindent\textsl{where formulae for the ideals }$\sigma_{p}(G)$\textsl{ has
been given by Theorem 3.5.}

\bigskip

\noindent\textbf{Proof.} It suffices to show that the inclusions
$\rho_{\left\vert \alpha\right\vert -1}\in H^{\ast}(G)$ induces an isomorphism
$\mathcal{F}(G)\cong\Delta_{\mathbb{Z}}(\rho_{\left\vert \alpha\right\vert
-1})_{\alpha\in S(G)}$.

By the proof of Theorem 2.5 the $\operatorname{mod}p$ reduction

\begin{quote}
$r_{p}:\mathbb{Z}[\omega_{i},y_{j}]_{1\leq i\leq n,1\leq j\leq m}%
\rightarrow\mathbb{F}_{p}[\omega_{i},y_{j}]_{_{1\leq i\leq n,t\in G(p)}}$
\end{quote}

\noindent satisfies the relations

\begin{quote}
$r_{p}(\alpha)=\left\{
\begin{tabular}
[c]{l}%
$\alpha^{(p)}$ if $\alpha=h_{i}$ or $g_{s}$ with $1\leq i\leq k$,
$s\in\overline{G}(p)$;\\
$\alpha_{j}^{(p)}\text{ (resp. }0\text{) if }\alpha=f_{j}$ with $j\in G(p)$
(resp. $j\in\overline{G}(p)$).
\end{tabular}
\ \right.  $
\end{quote}

\noindent It indicates that the reduction $r_{p}:E_{3}^{\ast,1}(G)\rightarrow
E_{3}^{\ast,1}(G;\mathbb{F}_{p})$ satisfies that

\begin{enumerate}
\item[(3.8)] $r_{p}(\rho_{\left\vert \alpha\right\vert -1})=\left\{
\begin{tabular}
[c]{l}%
$\xi_{\left\vert h_{i}^{(p)}\right\vert -1}$ if $\alpha=h_{i}$;\\
$-y_{j}^{k_{j}-1}\xi_{\left\vert \alpha_{j}^{(p)}\right\vert -1}$ (resp.
$p_{j}\xi_{\left\vert g_{j}^{(p)}\right\vert -1}$) if $\alpha=p_{j}g_{j}%
-y_{j}^{k_{j}-1}f_{j}$,
\end{tabular}
\ \ \ \right.  $
\end{enumerate}

\noindent where $j\in G(p)$ (resp. $j\in\overline{G}(p)$), and where the
classes $\xi_{\left\vert \alpha\right\vert -1}\in E_{3}^{\ast,1}%
(G;\mathbb{F}_{p})$ with $\alpha\in S_{p}(G)$ are defined in (3.1). By the
commutivity of the diagram induced by the monomorphism $\kappa$

\begin{quote}%
\begin{tabular}
[c]{lll}%
$E_{3}^{\ast,1}(G)$ & $\overset{r_{p}}{\rightarrow}$ & $E_{3}^{\ast
,1}(G;\mathbb{F}_{p})$\\
$\kappa\downarrow$ &  & $\kappa\downarrow$\\
$H^{\ast}(G)$ & $\overset{r_{p}}{\rightarrow}$ & $H^{\ast}(G;\mathbb{F}_{p})$%
\end{tabular}

\end{quote}

\noindent formula (3.8) holds for the action of $r_{p}$ on $H^{\ast}(G)$ as well.

Consider the monomial

\begin{quote}
$u=\prod\limits_{\alpha\in S(G)}\rho_{\left\vert \alpha\right\vert -1}\in$
$E_{3}^{g,n}(G)$($=H^{\dim G}(G)=\mathbb{Z}$ by (2.22))
\end{quote}

\noindent By (3.8) applying $r_{p}$ to $u$ yields the equality on $H^{\dim
G}(G;\mathbb{F}_{p})$

\begin{quote}
$r_{p}(u)\equiv(a\prod\limits_{s\in\overline{G}(p)}p_{s})\prod\limits_{t\in
G(p)}y_{t}^{k_{t}-1}\prod\limits_{\alpha\in S_{p}(G)}\xi_{\left\vert
\alpha\right\vert -1}\equiv(a\prod\limits_{s\in\overline{G}(p)}p_{s})u_{p}$,
\end{quote}

\noindent where $u_{p}\in E_{3}^{g,n}(G;\mathbb{F}_{p})(=H^{\dim
G}(G;\mathbb{F}_{p})=\mathbb{F}_{p}$ by (2.22)) is the class given in (3.3),
and where

\begin{quote}
$a=\left\{
\begin{tabular}
[c]{l}%
$(-1)^{\left\vert G(p)\right\vert }\text{, if either }G\neq E_{8}\text{ or
}G=E_{8}\text{, }p\neq2,5$;\\
$-1\text{, if }G=E_{8}\text{, }p=2$;\\
$2\text{, if }G=E_{8}\text{, }p=5$.
\end{tabular}
\ \ \ \ \ \ \right.  $
\end{quote}

\noindent Since the element $u_{p}$ generates the group $H^{\dim
G}(G;\mathbb{F}_{p})=\mathbb{F}_{p}$ by the proof of Lemma 3.3, and since the
coefficient $(a\prod_{s}p_{s})$ is always co--prime to $p$, the element
$r_{p}(u)$ generates $H^{\dim G}(G;\mathbb{F}_{p})$ for every prime $p$. It
follows that $u$ is a basis element of the group $E_{3}^{g,n}(G)=H^{\dim
G}(G)=\mathbb{Z}$. Since the class $u=\prod\limits_{\alpha\in S(G)}%
\rho_{\left\vert \alpha\right\vert -1}$ is also monotonous in the ring
$H^{\ast}(G)$ by Lemma 2.8, Lemma 2.9 concludes that the set of monomials

\begin{quote}
$\{\prod\rho_{\left\vert \alpha\right\vert -1}^{\varepsilon_{\alpha}}\in
H^{\ast}(G)\mid\alpha\in S(G)$, $\varepsilon_{\alpha}\in\{0,1\}\}$
\end{quote}

\noindent is linearly independent, and spans a direct summand of rank $2^{n}$
in a free part $\mathcal{F}(G)$ of $H^{\ast}(G)$. The proof is completed by
the relation $\dim H^{\ast}(G;\mathbb{R})=2^{n}$ by Theorem 3.4.$\square$

\bigskip

\noindent\textbf{Remark 3.7.} The proof of Theorem 3.6 indicates also the relation

\begin{quote}
$E_{3}^{\ast,\ast}(G)=\Lambda_{\mathbb{Z}}(\rho_{\left\vert \alpha\right\vert
-1})_{\alpha\in S(G)}\underset{p\in\{2,3,5\}}{\oplus}\operatorname{Im}%
\widehat{\beta}_{p}$.
\end{quote}

\noindent which implies that the same formula holds for $E_{\infty}^{\ast
,\ast}(G)$. Comparing this with (3.7) one obtains the additive isomorphisms

\begin{quote}
$E_{3}^{\ast,\ast}(G)=E_{\infty}^{\ast,\ast}(G)$, $E_{\infty}^{\ast,\ast
}(G)=H^{\ast}(G)$.
\end{quote}

\noindent conjectured by Ka\v{c} \cite{K1,K2} and Marlin \cite{M1}, respectively.

Results in Lemma 3.3 may be compared with \cite[Theorem 2, Theorem 3]{K2} by
Ka\v{c}. In our context the generators of the ring $E_{3}^{\ast,\ast
}(G;\mathbb{F}_{p})$ are fashioned explicitly from the presentation (1.4) of
the cohomology $H^{\ast}(G/T)$. As a result their degrees can be read directly
from the basic data of the corresponding Lie group $G$.$\square$

\section{The ring $H^{\ast}(G,\mathcal{R})$ of exceptional Lie groups}

Assume in this section that $G$ is a $1$--connected exceptional Lie group.
Based on Theorems 3.4 and 3.6 we recover the classical computation on the
rings $H^{\ast}(G;\mathbb{F})$ with $\mathbb{F}=\mathbb{R}$ or $\mathbb{F}%
_{p}$ in Section 4.1, and give a proof to Theorem 1.4 in Section 4.3.

Historically, the ring $H^{\ast}(G;\mathbb{F})$ with $\mathbb{F}=\mathbb{R}$
or $\mathbb{F}_{p}$ was calculated case by case, presented by generators with
different origins. As a result one could hardly analyzing the structure of the
integral cohomology $H^{\ast}(G)$ from information on $H^{\ast}(G;\mathbb{F}%
)$. In comparison, with our generators in various coefficients stemming solely
from the system $\{e_{i},f_{j},g_{j}\}$ in the Schubert classes on $G/T$, the
relationships between the cohomologies $H^{\ast}(G)$ and $H^{\ast
}(G;\mathbb{F}_{p})$ are transparent from the beginning. For this reason we
can proceed from Theorem 3.6 to decide the ring structure on $H^{\ast}(G)$.

In order to make the degrees of the generators transparent, the notion
$x_{\left\vert y_{i}\right\vert }$ is used instead of $\pi^{\ast}y_{i}$, where
$y_{1},\cdots,y_{m}$ are Schubert classes on $G/T$ given in Theorem 1.2 (see
Example 3.2). In addition, for a pair $(G;\mathcal{R})$ let $D(G;\mathcal{R})$
be the degree set of the primary polynomials of $G$ over $\mathcal{R}$,
arranged in the increasing order. We shall also set, for a prime
$p\in\{2,3,5\}$, that

\begin{enumerate}
\item[(4.1)] $d^{1}(G,p):=\{\left\vert \alpha_{t}^{(p)}\right\vert \mid$ $t\in
G(p)\}$, $\overline{d}^{1}(G,p):=\{\left\vert \beta_{t}^{(p)}\right\vert \mid$
$t\in G(p)\}$.
\end{enumerate}

\noindent In view of the inclusion (see Corollary 2.6 (resp. (3.6))

\begin{quote}
$d^{1}(G,p)\subset D(G;\mathbb{F}_{p})$ (resp. $\overline{d}^{1}(G,p)\subset
D(G;\mathbb{Z})$).
\end{quote}

\noindent we let $d^{2}(G,p)$ (resp. $\overline{d}^{2}(G,p)$) be the
complement of $d^{1}(G,p)$ in $D(G;\mathbb{F}_{p})$ (resp. $\overline{d}%
^{1}(G,p)$ in $D(G;\mathbb{Z})$).

\bigskip

\noindent\textbf{Example 4.1. }By the contents of Table 5.2 one finds that

1) for a pair $(G,p)$ with $G(p)\neq\emptyset$ the degree set $D(G,\mathbb{F}%
_{p})$ is given by the table below, where elements of the subset $d^{1}(G,p)$
are underlined:

\begin{center}
\noindent\textbf{Table 4.1.} The degree set $D(G,\mathbb{F}_{p})$ of primary
polynomials for $G(p)\neq\emptyset$

\qquad%
\begin{tabular}
[c]{l|l}\hline
$(G,p)$ & $d(G,p)\subset D(G,\mathbb{F}_{p})$\\\hline
$(G_{2},2)$ & $\{4,\underline{6}\}$\\
$(F_{4},2)$ & $\{4,\underline{6},16,24\}$\\
$(E_{6},2)$ & $\{4,\underline{6},10,16,18,24\}$\\
$(E_{7},2)$ & $\{4,\underline{6},\underline{10},16,\underline{18},24,28\}$\\
$(E_{8},2)$ & $\{4,\underline{6},\underline{10},16,\underline{18}%
,24,28,\underline{30}\}$\\
$(F_{4},3)$ & $\{4,\underline{8},12,16\}$\\
$(E_{6},3)$ & $\{4,\underline{8},10,12,16,18\}$\\
$(E_{7},3)$ & $\{4,\underline{8},12,16,20,28,36\}$\\
$(E_{8},3)$ & $\{4,\underline{8},16,\underline{20},28,36,40,48\}$\\
$(E_{8},5)$ & $\{4,\underline{12},16,24,28,36,40,48\}$\\\hline
\end{tabular}

\end{center}

2) in the remaining cases $\mathcal{R}=\mathbb{R}$, $\mathbb{Z}$ or
$\mathbb{F}_{p}$ with $G(p)=\emptyset$ one has $D(G,\mathcal{R})=$
$D(G,\mathbb{Z})$, where the set $D(G,\mathbb{Z})$ is presented in the table
below, and elements of the subset $\overline{d}^{1}(G,p)$ are indicated by the
subscript $(p)$:

\begin{center}
\textbf{Table 4.2.} The degree set $D(G,\mathbb{Z})$ of primary polynomials

\qquad%
\begin{tabular}
[c]{l|l}\hline\hline
Type of $G$ & $\ \ \overline{d}(G,p)\subset D(G,\mathbb{Z})$\\\hline
$G_{2}$ & $\{4,12_{(2)}\}$\\
$F_{4}$ & $\{4,12_{(2)},16,24_{(3)}\}$\\
$E_{6}$ & $\{4,10,12_{(2)},16,18,24_{(3)}\}$\\
$E_{7}$ & $\{4,12_{(2)},16,20_{(2)},24_{(3)},28,36_{(2)}\}$\\
$E_{8}$ & $\{4,16,24_{(3)},28,36_{(2)},40_{(2)},48_{(2)},60_{(2,3,5)}%
\}$\\\hline\hline
\end{tabular}
.$\square$
\end{center}

\subsection{The ring $H^{\ast}(G;\mathbb{F})$}

If $\mathbb{F}=\mathbb{R}$ or $\mathbb{F}_{p}$ with $G(p)=\emptyset$ one has
by Theorem 3.4 that

\begin{quote}
$H^{\ast}(G;\mathbb{F})=\Delta_{\mathbb{F}}(\xi_{s-1})_{s\in D(G,\mathbb{F})}$,
\end{quote}

\noindent where $p\neq2$ by Table 4.1. Since $u^{2}=0$ for $u\in
H^{odd}(G;\mathbb{F})$ with $\mathbb{F\neq F}_{2}$ the module $\Delta
_{\mathbb{F}}(\xi_{s-1})_{s\in D(G,\mathbb{F})}$ can be replaced by the ring
$\Lambda_{\mathbb{F}}(\xi_{s-1})_{s\in D(G,\mathbb{F})}$. Therefore, the
contents of Table 4.2 yields the following result that implies in particular
the classical computation of Yen \cite{Y}, Borel and Chevalley \cite{BC}

\bigskip

\noindent\textbf{Corollary 4.2.} \textsl{For} $\mathbb{F}=\mathbb{R}$
\textsl{or} $\mathbb{F}_{p}$ \textsl{with} $G(p)=\emptyset$ \textsl{we have}

\begin{quote}
$H^{\ast}(G_{2};\mathbb{F})=\Lambda_{\mathbb{F}}(\xi_{3},\xi_{11})$;$\qquad$

$H^{\ast}(F_{4};\mathbb{F})=\Lambda_{\mathbb{F}}(\xi_{3},\xi_{11},\xi_{15}%
,\xi_{23})$;

$H^{\ast}(E_{6};\mathbb{F})=\Lambda_{\mathbb{F}}(\xi_{3},\xi_{9},\xi_{11}%
,\xi_{15},\xi_{17},\xi_{23})$;

$H^{\ast}(E_{7};\mathbb{F})=\Lambda_{\mathbb{F}}(\xi_{3},\xi_{11},\xi_{15}%
,\xi_{19},\xi_{23},\xi_{27},\xi_{35})$;

$H^{\ast}(E_{8};\mathbb{F})=\Lambda_{\mathbb{F}}(\xi_{3},\xi_{15},\xi_{23}%
,\xi_{27},\xi_{35},\xi_{39},\xi_{47},\xi_{59})$.$\square$
\end{quote}

For an odd prime $p$ the elements of $H^{odd}(G;\mathbb{F}_{p})$ are also
square free. Therefore, combining Theorem 3.4 with the contents in Tables 4.1
yields the next results that imply the calculations \cite{A2,B2,B3} by Araki
and Borel.

\bigskip

\noindent\textbf{Corollary 4.3.} \textsl{For all the pairs }$(G,p)$\textsl{
with} $p\in\{3,5\}$ \textsl{and} $G(p)\neq\emptyset$\textsl{ we have}

\begin{quote}
$H^{\ast}(F_{4};\mathbb{F}_{3})=\mathbb{F}_{3}[x_{8}]/\left\langle x_{8}%
^{3}\right\rangle \otimes\Lambda_{\mathbb{F}_{3}}(\xi_{3},\xi_{7},\xi_{11}%
,\xi_{15})$;

$H^{\ast}(E_{6};\mathbb{F}_{3})=\mathbb{F}_{3}[x_{8}]/\left\langle x_{8}%
^{3}\right\rangle \otimes\Lambda_{\mathbb{F}_{3}}(\xi_{3},\xi_{7},\xi_{9}%
,\xi_{11},\xi_{15},\xi_{17})$;

$H^{\ast}(E_{7};\mathbb{F}_{3})=\mathbb{F}_{3}[x_{8}]/\left\langle x_{8}%
^{3}\right\rangle \otimes\Lambda_{\mathbb{F}_{3}}(\xi_{3},\xi_{7},\xi_{11}%
,\xi_{15},\xi_{19},\xi_{27},\xi_{35})$;

$H^{\ast}(E_{8};\mathbb{F}_{3})=\mathbb{F}_{3}[x_{8},x_{20}]/\left\langle
x_{8}^{3},x_{20}^{3}\right\rangle \otimes\Lambda_{\mathbb{F}_{3}}(\xi_{3}%
,\xi_{7},\xi_{15},\xi_{19},\xi_{27},\xi_{35},\xi_{39},\xi_{47})$;

$H^{\ast}(E_{8};\mathbb{F}_{5})=\mathbb{F}_{5}[x_{12}]/\left\langle x_{12}%
^{5}\right\rangle \otimes\Lambda(\xi_{3},\xi_{11},\xi_{15},\xi_{23},\xi
_{27},\xi_{35},\xi_{39},\xi_{47})$.$\square$
\end{quote}

\bigskip

For $p=2$ we have by Theorem 3.4 that

\begin{quote}
$H^{\ast}(G;\mathbb{F}_{2})=\mathbb{F}_{2}[x_{\left\vert y_{t}\right\vert
}]_{t\in G(2)}/\left\langle x_{\left\vert y_{t}\right\vert }^{k_{t}%
}\right\rangle \otimes\Delta_{\mathbb{F}_{2}}(\xi_{s-1})_{s\in D(G,\mathbb{F}%
_{2})}$.
\end{quote}

\noindent In view of this presentation the determination of the ring $H^{\ast
}(G;\mathbb{F}_{2})$ amounts to express all the squares $\xi_{s-1}^{2}$ with
$s\in D(G,\mathbb{F}_{2})$ as elements of the first factor $\mathbb{F}%
_{2}[x_{\left\vert y_{t}\right\vert }]_{t\in G(2)}/\left\langle x_{\left\vert
y_{t}\right\vert }^{k_{t}}\right\rangle $, see Lemma 2.8. For this purpose one
can make use of the Steenrod squares $Sq^{r}$, $r\geq1$, by which $\xi
_{s-1}^{2}=Sq^{s-1}\xi_{s-1}$ \cite{SE}.

With respect to the presentation (3.4) of the cohomology $H^{\ast
}(G;\mathbb{F}_{p})$ the structure of $H^{\ast}(G;\mathbb{F}_{p})$ as a module
over the Steenrod algebra has been determined in \cite{DZ3}. In particular, we
get from \cite[Corollary 4.4]{DZ3} that

\bigskip

\noindent\textbf{Corollary 4.4. }\textsl{For} $\mathbb{F}=\mathbb{F}_{2}$
\textsl{we have}

\begin{quote}
i) $H^{\ast}(G_{2};\mathbb{F}_{2})=\mathbb{F}_{2}[x_{6}]/\left\langle
x_{6}^{2}\right\rangle \otimes\Delta_{\mathbb{F}_{2}}(\xi_{3})\otimes
\Lambda_{\mathbb{F}_{2}}(\xi_{5})$;

ii) $H^{\ast}(F_{4};\mathbb{F}_{2})=\mathbb{F}_{2}[x_{6}]/\left\langle
x_{6}^{2}\right\rangle \otimes\Delta_{\mathbb{F}_{2}}(\xi_{3})\otimes
\Lambda_{\mathbb{F}_{2}}(\xi_{5},\xi_{15},\xi_{23})$;

iii) $H^{\ast}(E_{6};\mathbb{F}_{2})=\mathbb{F}_{2}[x_{6}]/\left\langle
x_{6}^{2}\right\rangle \otimes\Delta_{\mathbb{F}_{2}}(\xi_{3})\otimes
\Lambda_{\mathbb{F}_{2}}(\xi_{5},\xi_{9},\xi_{15},\xi_{17},\xi_{23})$;

iv) $H^{\ast}(E_{7};\mathbb{F}_{2})=\frac{\mathbb{F}_{2}[x_{6},x_{10},x_{18}%
]}{\left\langle x_{6}^{2},x_{10}^{2},x_{18}^{2}\right\rangle }\otimes
\Delta_{\mathbb{F}_{2}}(\xi_{3},\xi_{5},\xi_{9})\otimes\Lambda_{\mathbb{F}%
_{2}}(\xi_{15},\xi_{17},\xi_{23},\xi_{27})$;

v) $H^{\ast}(E_{8};\mathbb{F}_{2})=\frac{\mathbb{F}_{2}[x_{6},x_{10}%
,x_{18},x_{30}]}{\left\langle x_{6}^{8},x_{10}^{4},x_{18}^{2},x_{30}%
^{2}\right\rangle }\otimes\Delta_{\mathbb{F}_{2}}(\xi_{3},\xi_{5},\xi_{9}%
,\xi_{15},\xi_{23})\otimes\Lambda_{\mathbb{F}_{2}}(\xi_{17},\xi_{27},\xi
_{29})$,
\end{quote}

\noindent\textsl{where}

\begin{enumerate}
\item[(4.2)]
\begin{tabular}
[c]{l}%
$\xi_{3}^{2}=x_{6}$\textsl{ in }$G_{2},F_{4},E_{6},E_{7},E_{8}$\textsl{;}\\
$\xi_{5}^{2}=x_{10},\quad\xi_{9}^{2}=x_{18}$\textsl{ in }$E_{7},E_{8}%
$\textsl{;}\\
$\xi_{15}^{2}=x_{30}$\textsl{; }$\xi_{23}^{2}=x_{6}^{6}x_{10}$\textsl{ in
}$E_{8}$\textsl{.}$\square$%
\end{tabular}

\end{enumerate}

Historically, the ring $H^{\ast}(G;\mathbb{F}_{2})$ for the exceptional Lie
groups $G$ were obtained by Borel, Araki, Shikata and Kono
\cite{A1,AS,B2,B3,Ko} using "\textsl{transgressive generators}". Our
generators $\xi_{s-1}$ constructed in (3.1) may not be transgressive in the
sense of \cite{Ko}, compare Corollaries 4.2 and 4.4 in \cite{DZ3}.

\subsection{The ideal $\sigma_{p}(G)\subset H^{\ast}(G)$\ }

By Theorem 3.5 the ideal $\sigma_{p}(G)$ is non--trivial if and only if
$G(p)\neq\emptyset$. In addition, for a pair $(G,p)$ with $G(p)\neq\emptyset$
one has $d^{2}(G,p)=\overline{d}^{2}(G,p)$ by Example 4.1, while the formula
(3.8) implies that

\begin{enumerate}
\item[(4.3)] $r_{p}(\rho_{s-1})=\pm\xi_{s-1}$, $s\in d^{2}(G,p)$.
\end{enumerate}

\noindent On the other hand, if one defines for each subset $I\subset
d^{1}(G,p)$ the next elements of the integral cohomology $H^{\ast}(G)$ (in
analogue to (3.5))

\begin{enumerate}
\item[(4.4)] $\mathcal{C}_{I}=\beta_{p}(\xi_{I})$, $\mathcal{D}_{I}%
=\sum\limits_{t\in I}-x_{\left\vert y_{t}\right\vert }\mathcal{C}_{I_{t}}$,
$\mathcal{R}_{I}=(\prod\limits_{t\in I}^{k_{t}-1}x_{\left\vert y_{t}%
\right\vert }^{k_{t}-1})\mathcal{C}_{I_{t}}$,
\end{enumerate}

\noindent one has then the obvious relations

\begin{enumerate}
\item[(4.5)] $r_{p}(\mathcal{C}_{I})=c_{I}$; $r_{p}(\mathcal{D}_{I})=D_{I}$;
$r_{p}(\mathcal{R}_{I})=R_{I}$
\end{enumerate}

\noindent that allows us to translate the formula on $\operatorname{Im}%
\delta_{p}$ in ii) of Theorem 3.5 into the following additive presentation of
the ideal $\sigma_{p}(G)$

\begin{enumerate}
\item[(4.6)] $\sigma_{p}(G)=\frac{\mathbb{F}_{p}[x_{\left\vert y_{t}%
\right\vert }]_{t\in G(p)}/\left\langle x_{\left\vert y_{t}\right\vert
}^{k_{t}}\right\rangle \{1,\mathcal{C}_{I}\}^{+}}{\left\langle \mathcal{D}%
_{J},\mathcal{R}_{K}\right\rangle }\otimes\Delta(\rho_{s-1})_{s\in\overline
{d}^{2}(G,p)}$\textsl{, }
\end{enumerate}

\noindent where $I,J,K\subseteq d^{1}(G,p)$\ with $\left\vert I\right\vert
,\left\vert J\right\vert ,\left\vert K\right\vert \geq2$. In view of (4.6) the
determination of the ring structure on $\sigma_{p}(G)$ is reduced to

\begin{quote}
a) express the squares $\rho_{s-1}^{2}$ with $s\in\overline{d}^{2}(G,p)$ as
elements in the ring $\mathbb{F}_{p}[x_{\left\vert y_{t}\right\vert }]_{t\in
G(p)}/\left\langle x_{\left\vert y_{t}\right\vert }^{k_{t}}\right\rangle $,
see Lemma 2.8;

b) express the products $\mathcal{C}_{I}\cdot\mathcal{C}_{J}$ with
$I,J\subseteq d^{1}(G,p)$ as a $\frac{\mathbb{F}_{p}[x_{\left\vert
y_{t}\right\vert }]_{t\in G(p)}}{\left\langle x_{\left\vert y_{t}\right\vert
}^{k_{t}}\right\rangle }$--linear combinations of the $\mathcal{C}_{K}$'s.
\end{quote}

\noindent By identifying $\sigma_{p}(G)$ with the subring $\operatorname{Im}%
\delta_{p}\subset H^{\ast}(G;\mathbb{F}_{p})$ via $r_{p}$, the tasks in a) and
b) can be implemented by computation in the ring $H^{\ast}(G;\mathbb{F}_{p})$,
whose structure has already been settled by Corollaries 4.3 and 4.4. In the
following two results we note that $G(p)\neq\emptyset$ implies that
$p\in\{2,3,5\}$ by Theorem 1.2.

\bigskip

\noindent\textbf{Lemma 4.5.} \textsl{All the nontrivial torsion ideal }%
$\sigma_{p}(G)$\textsl{ with }$p=3,5$\textsl{ are}

\begin{quote}
i) $\sigma_{3}(F_{4})\cong\mathbb{F}_{3}[x_{8}]^{+}/\left\langle x_{8}%
^{3}\right\rangle \otimes\Lambda_{\mathbb{F}_{3}}(\varrho_{3},\varrho
_{11},\varrho_{15})$;

ii) $\sigma_{3}(E_{6})\cong\mathbb{F}_{3}[x_{8}]^{+}/\left\langle x_{8}%
^{3}\right\rangle \otimes\Lambda_{\mathbb{F}_{3}}(\varrho_{3},\varrho
_{9},\varrho_{11},\varrho_{15},\varrho_{17})$;

iii) $\sigma_{3}(E_{7})\cong\mathbb{F}_{3}[x_{8}]^{+}/\left\langle x_{8}%
^{3}\right\rangle \otimes\Lambda_{\mathbb{F}_{3}}(\varrho_{3},\varrho
_{11},\varrho_{15},\varrho_{19},\varrho_{27},\varrho_{35})$;

iv) $\sigma_{3}(E_{8})\cong\frac{\mathbb{F}_{3}[x_{8},x_{20},\mathcal{C}%
_{\{8,20\}}]^{+}}{\left\langle x_{8}^{3},x_{20}^{3},x_{8}^{2}x_{20}%
^{2}\mathcal{C}_{\{8,20\}},(\mathcal{C}_{\{8,20\}})^{2}\right\rangle }%
\otimes\Lambda_{\mathbb{F}_{3}}(\varrho_{3},\varrho_{15},\varrho_{27}%
,\varrho_{35},\varrho_{39},\varrho_{47})$;

v) $\sigma_{5}(E_{8})\cong\mathbb{F}_{5}[x_{12}]^{+}/\left\langle x_{12}%
^{5}\right\rangle \otimes\Lambda_{\mathbb{F}_{5}}(\varrho_{3},\varrho
_{15},\varrho_{23},\varrho_{27},\varrho_{35},\varrho_{39},\varrho_{47})$.
\end{quote}

\noindent\textbf{Proof.} With $p=3,5$ we have $u^{2}=0$ for $u\in
H^{odd}(G;\mathbb{F}_{p})$. It implies that the factor $\Delta(\rho
_{s-1})_{s\in\overline{d}^{2}(G,p)}$ in (4.6) can be replaced by the exterior
ring $\Lambda(\rho_{s-1})_{s\in\overline{d}^{2}(G,p)}$. One gets the
presentations i), ii), iii) and v) from (4.6), as well as the relations
$d^{1}(G,3)=\{8\}$ for $G=F_{4},E_{6},E_{7}$; $d^{1}(E_{8},5)=\{12\}$ by the
contents of Table 4.1. Note that the classes of the type $\mathcal{C}_{I}$ are
absent since in these cases the set $d^{1}(G,p)$ is always a singleton.

Similarly, one gets iv) from $d^{1}(E_{8},3)=\left\{  8,20\right\}  $, by
noticing that the element $\mathcal{C}_{\{8,20\}}$ is the only class of the
type $\mathcal{C}_{I}$ with $\left\vert I\right\vert \geq2$, whose square is
trivial for the degree reason.

Finally, the proof is completed by the relation $d^{1}(G_{2},3)=d^{1}%
(G,5)=\emptyset$, where $G\neq E_{8}$.$\square$

\bigskip

\noindent\textbf{Lemma 4.6.} \textsl{The rings }$\sigma_{2}(G)$ \textsl{for
the exceptional Lie groups }$G$\textsl{ are}

\begin{quote}
\textsl{i)} $\sigma_{2}(G_{2})=\mathbb{F}_{2}[x_{6}]^{+}/\left\langle
x_{6}^{2}\right\rangle \otimes\Delta_{\mathbb{F}_{2}}(\varrho_{3})$;

\textsl{ii)} $\sigma_{2}(F_{4})=\mathbb{F}_{2}[x_{6}]^{+}/\left\langle
x_{6}^{2}\right\rangle \otimes\Delta_{\mathbb{F}_{2}}(\varrho_{3}%
)\otimes\Lambda_{\mathbb{F}_{2}}(\varrho_{15},\varrho_{23})$;

\textsl{iii)} $\sigma_{2}(E_{6})=\mathbb{F}_{2}[x_{6}]^{+}/\left\langle
x_{6}^{2}\right\rangle \otimes\Delta_{\mathbb{F}_{2}}(\varrho_{3}%
)\otimes\Lambda_{\mathbb{F}_{2}}(\varrho_{9},\varrho_{15},\varrho_{17}%
,\varrho_{23})$;

\textsl{iv)} $\sigma_{2}(E_{7})=\frac{\mathbb{F}_{2}[x_{6},x_{10}%
,x_{18},\mathcal{C}_{I}]^{+}}{\left\langle x_{6}^{2},x_{10}^{2},x_{18}%
^{2},\mathcal{D}_{J},\mathcal{R}_{K},\mathcal{S}_{I,J}\right\rangle }%
\otimes\Delta_{\mathbb{F}_{2}}(\varrho_{3})\otimes\Lambda_{\mathbb{F}_{2}%
}(\varrho_{15},\varrho_{23},\varrho_{27})$ \textsl{with }$I,J,K\subseteq
d^{1}(E_{7},2)=\{6,10,18\}$\textsl{,} $\left\vert I\right\vert ,\left\vert
J\right\vert ,\left\vert K\right\vert \geq2$\textsl{;}

\textsl{v)} $\sigma_{2}(E_{8})=\frac{\mathbb{F}_{2}[x_{6},x_{10},x_{18}%
,x_{30},\mathcal{C}_{I}]^{+}}{\left\langle x_{6}^{8},x_{10}^{4},x_{18}%
^{2},x_{30}^{2},\mathcal{D}_{J},\mathcal{R}_{K},\mathcal{S}_{I,J}\right\rangle
}\otimes\Delta_{\mathbb{F}_{2}}(\varrho_{3},\varrho_{15},\varrho_{23}%
)\otimes\Lambda_{\mathbb{F}_{2}}(\varrho_{27})$ \textsl{with }$I$\textsl{,}
$J,K\subseteq d^{1}(E_{8},2)=\{6,10,18,30\}$\textsl{,} $\left\vert
I\right\vert ,\left\vert J\right\vert ,\left\vert K\right\vert \geq
2$\textsl{,}
\end{quote}

\noindent\textsl{where }

\begin{quote}
$\varrho_{3}^{2}=x_{6}$ \textsl{for all} $G$\textsl{, }$\varrho_{15}%
^{2}=x_{30},\varrho_{23}^{2}=x_{6}^{6}x_{10}$ \textsl{for} $G=E_{8}$\textsl{,}
\end{quote}

\noindent\textsl{and} \textsl{where the relations of the type }$\mathcal{S}%
_{I,J}$\textsl{ in iv) and v) is}

\begin{enumerate}
\item[(4.7)] $\mathcal{S}_{I,J}=\mathcal{C}_{I}\mathcal{C}_{J}+\sum
\limits_{t\in I}x_{t}\prod\limits_{s\in I_{t}\cap J}\xi_{s-1}^{2}%
\mathcal{C}_{\left\langle I_{t},J\right\rangle }$
\end{enumerate}

\noindent\textsl{with} $\left\langle I,J\right\rangle =\{t\in I\cup J\mid
t\notin I\cap J\}$\textsl{, }$\prod_{s\in I_{t}\cap J}\xi_{s-1}^{2}=1$
\textsl{if }$I_{t}\cap J=\emptyset$\textsl{, and with} \textsl{the squares
}$\xi_{s-1}^{2}$\textsl{ being evaluated by formulae (4.2).}

\bigskip

\noindent\textbf{Proof. }We can focus on the relatively nontrivial case
$G=E_{8}$ for which the formula (4.6) turns to be

\begin{quote}
a) $\sigma_{2}(E_{8})=\frac{\mathbb{F}_{2}[x_{6},x_{10},x_{18},x_{30}%
]\{1,\mathcal{C}_{I}\}^{+}}{\left\langle x_{6}^{8},x_{10}^{4},x_{18}%
^{2},x_{30}^{2},\mathcal{D}_{J},\mathcal{R}_{K}\right\rangle }\otimes
\Delta_{\mathbb{F}_{2}}(\varrho_{3},\varrho_{15},\varrho_{23},\varrho_{27})$,
\end{quote}

\noindent where $I,J,K\subseteq d^{1}(E_{8},2)=\{6,10,18,30\}$\textsl{,}
$\left\vert I\right\vert ,\left\vert J\right\vert ,\left\vert K\right\vert
\geq2$. Since the reduction $r_{2}$ restricts to an isomorphism

\begin{quote}
$\sigma_{2}(E_{8})\cong\operatorname{Im}\delta_{2}\subset H^{\ast}%
(E_{8};\mathbb{F}_{2})$ (by i) of Theorem 3.5)
\end{quote}

\noindent one has by (4.3) and (4.2) that

\begin{quote}
$\Delta_{\mathbb{F}_{2}}(\varrho_{3},\varrho_{15},\varrho_{23},\varrho
_{27})=\Delta_{\mathbb{F}_{2}}(\varrho_{3},\varrho_{15},\varrho_{23}%
)\otimes\Lambda_{\mathbb{F}_{2}}(\varrho_{27})$,

\end{quote}

\noindent and that $\varrho_{3}^{2}=x_{6},\varrho_{15}^{2}=x_{30},\varrho
_{23}^{2}=x_{6}^{6}x_{10}$.

It remains to decide the multiplicative rule (4.7) among the classes
$\mathcal{C}_{I}$'s. For $I,J\subseteq d^{1}(E_{8},2)$ with $\left\vert
I\right\vert ,\left\vert J\right\vert \geq2$ we have in the ring $H^{\ast
}(E_{8};\mathbb{F}_{2})$ that

\begin{quote}
b) $\delta_{2}(\xi_{I})=\sum\limits_{t\in I}x_{t}\xi_{I_{t}}$ (by relation a)
in the proof of Theorem 3.5);

c) $\xi_{I}\xi_{J}=\prod\limits_{s\in I\cap J}\xi_{s-1}^{2}\xi_{\left\langle
I,J\right\rangle }$ with $\prod\limits_{s\in I\cap J}\xi_{s-1}^{2}=1$ if
$I\cap J=\emptyset$.
\end{quote}

\noindent Granted with i) of Theorem 3.5 the relation (4.7) is shown by

\begin{quote}
$c_{I}=\delta_{2}(\xi_{I})$ (see (3.5)), $r_{2}(\mathcal{C}_{I})=c_{I}$ (see (4.5)),
\end{quote}

\noindent as well as the calculation in the ring $H^{\ast}(E_{8}%
;\mathbb{F}_{2})$ (see v) of Corollary 4.4):

\begin{quote}
$c_{I}c_{J}=\delta_{2}(\xi_{I})\delta_{2}(\xi_{J})=\delta_{2}(\delta_{2}%
(\xi_{I})\xi_{J})$ (since $\delta_{2}^{2}=0$)

$=\delta_{2}(\sum\limits_{t\in I}x_{t}\xi_{I_{t}}\xi_{J})$ (by b))

$=\delta_{2}(\sum\limits_{t\in I}x_{t}\prod\limits_{s\in I_{t}\cap J}\xi
_{s-1}^{2}\xi_{\left\langle I_{t},J\right\rangle })$ (by c))

$=\sum\limits_{t\in I}x_{t}\prod\limits_{s\in I_{t}\cap J}\xi_{s-1}%
^{2}c_{\left\langle I_{t},J\right\rangle }$ (by $\delta_{2}(x_{t})=\delta
_{2}(\xi_{s-1}^{2})=0$)
\end{quote}

\noindent This establishes the presentation v) of the ring $\sigma_{2}(E_{8}%
)$.$\square$

\bigskip

\noindent\textbf{Example 4.7.} Taking $G=E_{8}$ as an example and noticing
that $d^{1}(E_{8},2)=\{6,10,18,30\}$, we have by (4.7) that

\begin{quote}
a) $\mathcal{C}_{\{6,10\}}\mathcal{C}_{\{6,10\}}=x_{6}\xi_{9}^{2}%
\mathcal{C}_{\{6\}}+x_{10}\xi_{5}^{2}\mathcal{C}_{\{10\}}=x_{6}^{2}%
x_{18}+x_{10}^{3}$
\end{quote}

\noindent since $(\xi_{9}^{2},\xi_{5}^{2})=(x_{18},x_{10})$, $(c_{\{6\}}%
,c_{\{10\}})=(x_{6},x_{10})$ in $H^{\ast}(E_{8};\mathbb{F}_{2})$;

\begin{quote}
b) $\mathcal{C}_{(6,10)}\mathcal{C}_{(6,18)}=x_{10}\xi_{5}^{2}\mathcal{C}%
_{\{18\}}+x_{6}\mathcal{C}_{\{6,10,18\}}=x_{10}^{2}x_{18}+x_{6}\mathcal{C}%
_{\{6,10,18\}}$
\end{quote}

\noindent since $\xi_{5}^{2}=x_{10}$ and $c_{\{18\}}=x_{18}$ in $H^{\ast
}(E_{8};\mathbb{F}_{2})$. These computation indicate that the multiplicative
rule (4.7) on $\sigma_{2}(E_{8})$ is highly non--trivial, though can be easily
implemented.$\square$

\subsection{The ring $H^{\ast}(G)$ (the proof of Theorem 1.4)}

By formula (3.7) the integral cohomology $H^{\ast}(G)$ has the following
presentation in which the subrings $\sigma_{p}(G)$ have been determined in
Section 4.2

\begin{enumerate}
\item[(4.8)] $H^{\ast}(G)=\Delta_{\mathbb{Z}}(\varrho_{s-1})_{s\in
D(G,\mathbb{Z})}\underset{p\in\{2,3,5\}}{\oplus}\sigma_{p}(G)$,
\end{enumerate}

\noindent Moreover, since $\sigma_{p}(G)$ is an ideal the product on $H^{\ast
}(G)$ defines an action

\begin{quote}
$\Delta_{\mathbb{Z}}(\varrho_{s-1})_{s\in D(G,\mathbb{Z})}\otimes\sigma
_{p}(G)\rightarrow\sigma_{p}(G)$
\end{quote}

\noindent of the free part on $\sigma_{p}(G)$. Therefore, to clarify the ring
structure on $H^{\ast}(G)$ it remains to decide the formulae that express

\begin{quote}
1) the squares $\varrho_{s-1}^{2}$ with $s\in D(G;\mathbb{Z})$ as elements of
the subring $\pi^{\ast}E_{3}^{\ast,0}(G,\mathbb{Z})\subset H^{\ast}(G)$ (see
Lemma 2.8, Example 3.2);

2) all the products $\varrho_{s-1}\cdot\mathcal{C}_{I}$ with $s\in\overline
{d}^{1}(G;p)$ and $I\subseteq d^{1}(G;p)$ as elements of the ideal $\sigma
_{p}(G)$.
\end{quote}

\noindent The relations on $H^{\ast}(G)$ that implement these two requests are
denoted by $\mathcal{F}_{s}$ and $\mathcal{H}_{s,I}$, respectively, and are
made explicit in the following result.

In term of (4.2) one has the one to one correspondence

\begin{quote}
$\varepsilon:\overline{d}^{1}(G;p)\rightarrow d^{1}(G;p)$ by $\varepsilon
(\left\vert \beta_{t}^{(p)}\right\vert )=\left\vert \alpha_{t}^{(p)}%
\right\vert $, $t\in G(p)$.
\end{quote}

\noindent Alternatively, the map $\varepsilon$ is characterized uniquely by
the relation

\begin{quote}
$\left\vert \beta_{t}^{(p)}\right\vert =k_{t}\varepsilon(\left\vert \beta
_{t}^{(p)}\right\vert )$, $t\in G(p)$,
\end{quote}

\noindent see Theorem 1.2.

\bigskip

\noindent\textbf{Lemma 4.8.} \textsl{The relations of the type }%
$\mathcal{F}_{s}$\textsl{ with }$s\in D(G;\mathbb{Z})$ \textsl{are given by
}$\varrho_{s-1}^{2}=0$\textsl{ with the following three exceptions:}

\begin{quote}
$\varrho_{3}^{2}=x_{6}$\textsl{ for all }$G$\textsl{;}

$\varrho_{15}^{2}=x_{30},$\textsl{ }$\varrho_{23}^{2}=x_{6}^{6}x_{10}$\textsl{
for }$G=E_{8}$\textsl{.}
\end{quote}

\textsl{The relations of the type }$\mathcal{H}_{s,I}$ \textsl{with}
$s\in\overline{d}^{1}(G;p)$\textsl{,} $I\subseteq d^{1}(G;p)$\textsl{ are
given by the following three possibilities}

\begin{enumerate}
\item[(4.9)] $\varrho_{s-1}\cdot\mathcal{C}_{I}=\left\{
\begin{tabular}
[c]{l}%
$x_{\varepsilon(s)}^{\frac{s}{\varepsilon(s)}-1}\cdot\mathcal{C}%
_{I\cup\{\varepsilon(s)\}}$\textsl{,} \textsl{if} $\varepsilon(s)\notin
I$\textsl{;}\\
$0$\textsl{,} \textsl{if either} $\varepsilon(s)\in I$\textsl{,} $p$
\textsl{is odd or} $I=\{\varepsilon(s)\}$, $p=2$\textsl{;}\\
$x_{\varepsilon(s)}^{\frac{s}{\varepsilon(s)}-1}\cdot\xi_{\varepsilon
(s)-1}^{2}\cdot\mathcal{C}_{I_{\varepsilon(s)}}$\textsl{,} \textsl{if} $p=2$,
$\varepsilon(s)\in I$ \textsl{and} $\left\vert I\right\vert \geq2$\textsl{,}%
\end{tabular}
\ \ \ \ \ \ \ \ \ \ \ \ \ \ \right.  $
\end{enumerate}

\noindent\textsl{where in the third instance the squares }$\xi_{\varepsilon
(s)-1}^{2}$\textsl{ have been evaluated by (4.2).}

\bigskip

\noindent\textbf{Proof.} According to (3.9) we have

\begin{center}
$r_{2}(\varrho_{s-1}^{2})\equiv\left\{
\begin{tabular}
[c]{l}%
$\xi_{s-1}^{2}$, if $s\notin\overline{d}^{1}(G,q)$ for all $q\in\{2,3,5\}$;\\
$q^{2}\cdot\xi_{s-1}^{2}$, if $s\in\overline{d}^{1}(G,q)$ for some $q\neq2$;\\
$x_{\varepsilon(s)}^{2(\frac{s}{\varepsilon(s)}-1)}\cdot\xi_{\varepsilon
(s)-1}^{2}=0$ if $s\in\overline{d}^{1}(G,2)$ (since $x_{\varepsilon(s)}%
^{\frac{s}{\varepsilon(s)}}\equiv0$)
\end{tabular}
\ \ \ \ \right.  $.
\end{center}

\noindent The relations of the type $\mathcal{F}_{s}$ are verified by

a) $\varrho_{s-1}^{2}\in\sigma_{2}(G)$;

b) $r_{2}$ restricts to an isomorphism $\sigma_{2}(G)\rightarrow
\operatorname{Im}\delta_{2}\subset H^{\ast}(G;\mathbb{F}_{2})$;

c) the results on $\xi_{s-1}^{2}$ given in (4.2).

Similarly, resorting to the isomorphism $r_{p}:\sigma_{p}(G)\cong
\operatorname{Im}\delta_{p}$ by i) of Theorem 3.5 the relation (4.9) of the
type $\mathcal{H}_{s,I}$ with $s\in\overline{d}^{1}(G;p)$ and $I\subseteq
d^{1}(G;p)$ are obtained by the calculation

\begin{quote}
$r_{p}(\varrho_{s-1}\mathcal{C}_{I})\equiv-x_{\varepsilon(s)}^{\frac
{s}{\varepsilon(s)}-1}\cdot\xi_{\varepsilon(s)-1}c_{I}$ (by (3.9) and
$r_{p}(\mathcal{C}_{I})=c_{I}$)

$\equiv-x_{\varepsilon(s)}^{\frac{s}{\varepsilon(s)}-1}\cdot\xi_{\varepsilon
(s)-1}\delta_{p}(\xi_{I})$ (by $c_{I}=\delta_{p}(\xi_{I})$)

$\equiv-x_{\varepsilon(s)}^{\frac{s}{\varepsilon(s)}-1}\delta_{p}%
(\xi_{\{\varepsilon(s)\}}\xi_{I})$ (since $x_{\varepsilon(s)}^{\frac
{s}{\varepsilon(s)}-1}\delta_{p}(\xi_{\{\varepsilon(s)\}})\equiv
x_{\varepsilon(s)}^{\frac{s}{\varepsilon(s)}}\equiv0$)

$\equiv\left\{
\begin{tabular}
[c]{l}%
$-x_{\varepsilon(s)}^{\frac{s}{\varepsilon(s)}-1}c_{I\cup\{\varepsilon
(s)\}}\text{, if }\varepsilon(s)\notin I\text{;\quad}$\\
$0\text{, if }I=\{\varepsilon(s)\}\text{;}$\\
$-x_{\varepsilon(s)}^{\frac{s}{\varepsilon(s)}-1}(\xi_{\varepsilon(s)-1}%
)^{2}c_{I_{\varepsilon(s)}}\text{, if }\varepsilon(s)\in I\text{, }\left\vert
I\right\vert \geq2$.
\end{tabular}
\ \ \ \ \right.  $
\end{quote}

\noindent Note that in the third instance we must have $(\xi_{\varepsilon
(s)-1})^{2}=0$ when $p\neq2$ since $\xi_{\varepsilon(s)-1}$ is of odd degree
with order $\neq2$, and that the squares $(\xi_{\varepsilon(s)-1})^{2}$ when
$p=2$ have been evaluated by (4.2).$\square$

\bigskip

We are ready to show Theorem 1.4 stated in Section 1.

\bigskip

\noindent\textbf{Proof of Theorem 1.4.} With respect to the presentation (4.8)
the structure of $H^{\ast}(G)$ as a ring has been decided by Lemmas 4.5, 4.6
and 4.8. Concerning the presentations (1.8)--(1.12) in Theorem 1.4 we need
only to remark that

i) an element $\varrho_{s-1}$ with free square contributes to a generator in
the exterior factor of the free part;

ii) if the set $d^{1}(G;p)$ is a singleton, the relation of the type
$\mathcal{H}_{s,I}$ with $s\in\overline{d}^{1}(G;p)$, $I\subseteq d^{1}(G;p)$
is unique, and can be presented as $x_{\varepsilon(s)}\varrho_{s-1}%
=0$.$\square$

\section{The basic data of simple Lie groups}

In the course of establishing Theorem 1.2 in \cite{DZ4}\textbf{ }a set
$\left\{  y_{1},\cdots,y_{m}\right\}  $ of special Schubert classes, as well
as the corresponding system $\left\{  h_{i},f_{j},g_{j}\right\}  _{1\leq i\leq
k,1\leq j\leq m}$ of polynomials, have been made explicit for every simple Lie
group, see \cite[Section 6]{DZ4}. As a result we obtains the basic data
\cite[Tables 6.1, 6.2]{DZ4} of the simple Lie groups that are recorded below,
by which one obtains the degree sets $D(G;\mathcal{R})$ of the primary
polynomials of $G$ over $\mathcal{R}$ in Example 4.1.

\bigskip

\noindent\textbf{Table 5.1.} Basic data for the classical groups

\begin{center}
{\footnotesize
\begin{tabular}
[c]{l|llll}\hline\hline
$G$ & $SU(n+1)$ & $Sp(n)$ & $Spin(2n)$ & $Spin(2n+1)$\\\hline
$\{k,m\}$ & $\{n,0\}$ & $\{n,0\}$ & $\{[\frac{n+3}{2}],[\frac{n-2}{2}]\}$ &
$\{[\frac{n+2}{2}],[\frac{n-1}{2}]\}$\\
$\{\deg h_{i}\}$ & $\{2i+2\}$ & $\{4i\}$ & $\{4t,2n,2^{[\log_{2}%
(n-1)]+2}\}_{1\leq t\leq\lbrack\frac{n-1}{2}]}$ & $\{4t,2^{[\log_{2}%
n]+2}\}_{1\leq t\leq\lbrack\frac{n}{2}]}$\\
$\{\deg y_{j}\}$ &  &  & $\{4j+2\}$ & $\{4j+2\}$\\
$\{p_{j}\}$ &  &  & $\{2,\cdots,2\}$ & $\{2,\cdots,2\}$\\
$\{k_{j}\}$ &  &  & $\{2^{[\log_{2}\frac{n-1}{2j+1}]+1}\}$ & $\{2^{[\log
_{2}\frac{n}{2j+1}]+1}\}$\\\hline\hline
\end{tabular}
}
\end{center}

\noindent\textbf{Table 5.2.} Basic data for the exceptional Lie groups.

\begin{center}
{\footnotesize
\begin{tabular}
[c]{l}\hline\hline%
\begin{tabular}
[c]{l|lllll}%
$G$ & $G_{2}$ & $F_{4}$ & $E_{6}$ & $E_{7}$ & $E_{8}$\\\hline
$\{k,m\}$ & $\{1,1\}$ & $\{2,2\}$ & $\{4,2\}$ & $\{3,4\}$ & $\{3,7\}$\\
$\{\deg h_{i}\}$ & $\{4\}$ & $\{4,16\}$ & $\{4,10,16,18\}$ & $\{4,16,28\}$ &
$\{4,16,28\}$\\
$\{\deg y_{j}\}$ & $\{6\}$ & $\{6,8\}$ & $\{6,8\}$ & $\{6,8,10,18\}$ &
$\{6,8,10,12,18,20,30\}$\\
$\{p_{j}\}$ & $\{2\}$ & $\{2,3\}$ & $\{2,3\}$ & $\{2,3,2,2\}$ &
$\{2,3,2,5,2,3,2\}$\\
$\{k_{j}\}$ & $\{2\}$ & $\{2,3\}$ & $\{2,3\}$ & $\{2,3,2,2\}$ &
$\{8,3,4,5,2,3,2\}$%
\end{tabular}
\\\hline\hline
\end{tabular}
}
\end{center}

\end{document}